\documentclass[11pt,reqno,twoside]{article}
\usepackage{amsmath,amsfonts,amssymb}
\usepackage{mathrsfs}
\usepackage{amsthm}
\usepackage{enumitem}
\usepackage[utf8]{inputenc}
\usepackage{fontenc}
\usepackage{hyperref}

 \usepackage{subfig}
\usepackage[english]{babel}
\usepackage{amsbsy,amscd,fancyhdr,graphicx,psfrag,fancybox,indentfirst,color}
\usepackage{graphics,epsfig}
\usepackage{mysty}
\usepackage{fullpage}
\usepackage{graphicx}

\title{\textbf {Nonlocal $p$-Laplacian evolution problems on graphs}}         
\date{}
\author{Yosra Hafiene, Jalal Fadili and Abderrahim Elmoataz \\\\
\small{Normandie Univ, ENSICAEN, UNICAEN, CNRS, GREYC, France}}

\begin{document}

\maketitle

\begin{abstract}
In this paper we study numerical approximations of the evolution problem for the nonlocal $p$-Laplacian with homogeneous Neumann boundary conditions. First, we derive a bound on the distance between two continuous-in-time trajectories defined by two different evolution systems (i.e. with different kernels and initial data). We then provide a similar bound for the case when one of the trajectories is discrete-in-time and the other is continuous. In turn, these results allow us to establish error estimates of the discretized $p$-Laplacian problem on graphs. More precisely, for networks on convergent graph sequences (simple and weighted graphs), we prove convergence and provide rate of convergence of solutions for the discrete models to the solution of the continuous problem as the number of vertices grows. We finally touch on the limit as $p \to \infty$ in these approximations and get uniform convergence results.

\end{abstract}




\section{Introduction}
\subsection{Problem formulation}

Our main goal in this paper is to study the following nonlinear diffusion problem, which we call the nonlocal $p$-Laplacian problem with homogeneous Neumann boundary conditions: 
 \begin{equation}\tag{\textrm{$\mathcal{P}$}}
  \begin{cases}
  u_t(x,t) = \frac{\partial }{\partial t} u(x,t) = -\A(u(x,t)), & x \in \O, t > 0,
  \\
  u(x,0)= g(x), \quad x \in \O,
  \end{cases}
  \label{neumann}
  \end{equation}
where \newcor{$ p \in ]1, + \infty[$} and 
\[
\A(u(x,t)) =  - \int_{\O} \K(x,y) \abs{u(y,t) - u(x,t) }^{p-2} (u(y,t) - u(x,t) ) dy,
\]
$\O \subset \mathbb{R} $ is a bounded domain, without loss of generality $\O = [0,1]$, and $\K(\cdot, \cdot)$ is a symmetric, nonnegative and bounded \newcor{function}. In particular, the kernel $\K(\cdot,\cdot)$ represents the limit object for some convergent graph sequence $\{G_n\}_n, n \in \N$ for every $(x,y) \in \O^2$, whose meaning and form will be specified in the sequel, separately for every class of problems that we consider below. 

\newcor{The chief goal of this paper is to study numerical approximations of the evolution problem \eqref{neumann}, which in turn, will allow us to establish consistency estimates of the fully discretized $p$-Laplacian problem on graphs.} In recent years, partial differential equations (PDEs) involving the nonlocal $p$-Laplacian operator have become more and more \newcor{popular} both in the setting of Euclidean domains and on discrete graphs, as the $p$-Laplacian problem has been possessing many important features shared by many practical problems in mathematics, physics, engineering, biology, and economy, such as continuum mechanics, phase transition phenomena, population dynamics, see~\cite{mazon, mechanics} and references therein. Some closely related applications can be found in image processing, such as spectral clustering~\cite{cluster}, computer vision and machine learning~\cite{image, international, elliptic}.

Particularly, if $\K(x,y) = J(x-y)$, where the kernel $J : \R^d \to \R$ \newcor{is a nonnegative continuous radial function with compact support verifying $J(0) > 0$ and $\int_{\R^d} J(x) dx = 1$}, nonlocal evolution equations of the form 
\begin{equation*}
u_t(x,t) = J \ast u (x,t) - u(x,t) = \int_{\R^d} J(x-y) (u(y,t) - u(x,t)) dy,
\end{equation*}
where $\ast$ stands for the convolution, have many applications in modeling diffusion processes. See, among many others references,~\cite{mazon, phase, waves, effect, like, patterns, interfacial}.
As stated in~\cite{like},  in modeling the dispersal of organisms in space when $u(x,t)$ is their density at the point $x$ at time $t$, $J(x-y)$ is considered as the probability distribution of jumping from position $y$ to position $x$, then, the expression $J *u - u$ represents transport due to long-range dispersal mechanisms, that is the rate at which organisms are arriving to location $x$ from any other place.

Let us note that, with the definition of the solution, the evolution problem~\eqref{neumann} is the gradient flow associated to the functional 
\begin{equation}
\label{eq:Fp}
F_p(v)= \frac{1}{2p} \int_{\O^2} \K(x,y)\abs{v(y) - v(x)}^p dy dx,
\end{equation}
which is the nonlocal analog to the energy functional $ \int_{\O} \abs{   \nabla v}^p$ associated to the local $p$-Laplacian.

%
%

Solutions of \eqref{neumann} will  be understood in the following sense: 
\begin{defi}\label{def:sol}
A solution of~\eqref{neumann} in $[0,T]$ is a function 
\[
u \in W^{1,1} (0,T;L^1(\O)),
\]
that satisfies $u(x,0) = g(x)$ a.e. $x \in \O$ and 
\[
 u_t(x,t) = -\A(u(x,t)) \quad \text{ a.e. in } \O \times ]0,T[ . 
 \]
\end{defi}

\begin{rem}
\label{rem:sol}
 Observe that since $u \in W^{1,1} (0,T;L^1(\O))$, we have that $u$ is also a \textbf{strong} solution (see \cite[Definition~A.3]{vaillo}). Indeed,
 \begin{equation*}
 \vspace{1cm}
 \left.\begin{matrix}
C(0,T;L^1(\O)) \subset W^{1,1} (0,T;L^1(\O)) \vspace{0.5cm}\\  
W^{1,1} (0,T;L^1(\O)) \subset W_{loc}^{1,1} (0,T;L^1(\O))
\end{matrix}\right\} \Rightarrow  u \in C(0,T;L^1(\O)) \cap W_{loc}^{1,1} (0,T;L^1(\O)).
 \end{equation*}
\end{rem}

\subsection{Contributions}
\newcor{In this work we \newcor{intend} to provide two related contributions. Their combination, associated with techniques from the recent theory of graph limits, allow to quantitatively analyze evolution problems on convergent graph sequences and their limiting behaviour.}

More precisely, we first study the convergence and stability properties of the numerical solutions for the general time-continuous problem valid uniformly for $t \in [0,T]$, where $T> 0$. Under the assumption $p \in ]1, + \infty [$, as $n \to \infty$, we prove that the solution to this problem, that can be regarded as a spatial semi-discrete approximation of the initial problem via the kernel discretization, converges to a nonlocal evolution problem. We give Kobayashi-type estimates. Then, we apply our analysis to the forward and backward Euler schemes to get similar estimates for the fully discretized problem. In addition, we obtain convergence in the $L^p(\O)$ norm for both time continuous and totally  discretized problems. Convergence in $L^2(\O)$ norm is thus a corollary. We obtain these results without any extra regularity assumption.

Secondly, we apply these results to dynamical networks on simple and weighted graphs to show that the approximation of solutions of the discrete problems on simple and weighted graph sequences converge to those of the continuous problem. We give also a rate of convergence estimate. Specifically, for simple graph sequences, we show how the accuracy of the approximation depends on the regularity of the boundary of support of the graph limit \newcor{in the same vein as~\cite{medv} who did it for a nonlocal nonlinear heat equation (see also discussion in the forthcoming section)}. \newcor{In addition, for weighted graphs, we give a precise error estimate under the mild assumption that both the kernel $K$ and the initial data $g$ are also in Lipschitz spaces, which in particular contain functions of bounded variation.}

Let us note that we look in detail to the one-dimensional case, that is $\O = [0,1]$,  our results also hold when we deal with approximations in a multidimensional domain, since the extension to larger dimension spaces is straightforward. The proofs are similar to the one-dimensional case and are left to the reader.

\subsection{Relation to prior work}
Concerning previous work for this model, \newcor{the authors of~\cite{rossi}} have already obtained a similar conclusion under different but complementary assumptions. \newcor{Indeed, in that paper, only the case $\K(x,y) = J(x-y)$ was considered}. \newcor{The authors} showed that solutions to the numerical scheme converge to the continuous solution for both semi-discrete and totally discrete approximations. However, the convergence is only uniform and requires the positivity of the solution. 

\newcor{Another closely related and important work is that in~\cite{medv,medvedevrandom} which paved the way to study limit phenomena of evolution problems on both deterministic and random graphs. In~\cite{medv}, the author focused on a nonlinear heat equation on graphs, where the function $\F$ (see the proof of Theorem~\ref{main}) was assumed Lipschitz-continuous. This assumption was essential to prove well-posedness (existence and uniqueness follow immediately from the contraction principle), as well as to study the consistency in $L^2$ of the spatial semi-discrete approximation on simple and weighted graph sequences. Though this seminal work was quite inspiring to us, it differs from our work in many crucial aspects. First, the nonlocal $p$-Laplacian evolution problem at hand is different and cannot be covered by~\cite{medv} where $\F$ lacks Lipschitzianity, and thus raises several challenges (including for well-posedness and error estimates). Our results on Kobayashi-type estimates are also novel and are of independent interest beyond problems on networks. We also consider both the semi-discrete and fully-discrete versions with both forward and backward Euler approximations, that we fully characterize.}

\subsection{Paper organization}
This paper is organized as follows. In Section~\ref{graphlimitss}, we start with a general review of the necessary background on graph limits and represent the different types of graphs that we are going to deal with later. In Section~\ref{existenceuniqueness}, we address the well-posedness of the problem~\eqref{neumann}, we show that~\eqref{neumann} admits a unique solution in $C(0,T;L^1(\O))$. Further, in Theorem~\ref{existence}  we give a steadiness condition regarding the stability of the solution with respect to the initial data, which guarantees that the solution of~\eqref{neumann} remains in $L^p(\O)$, $1<p< +\infty$ as long as the initial data is in this same space. In particular, we apply this result to get our estimate bounds in the subsequent sections. In Section~\ref{timecontinuous} and~\ref{timediscrete}, we study the consistency of the time-continuous and time-discrete problems, respectively,  and establish some error estimates. Here, we extend~\eqref{neumann} to get the problem~\eqref{general} that we keep in mind as a space-discretized version of~\eqref{neumann} via the discretization of the kernel $K$, since we have the idea of applying it to study the relation between the solutions of the totally discrete problems~\eqref{neumanndisc} and~\eqref{pond} corresponding to simple and weighted graph sequences, respectively, and that of the initial problem~\eqref{neumann}, which is the subject of section~\ref{applications}. In Section~\ref{simplegraphs}, for sequences of simple graphs converging to $\{0,1\}$-valued graphons, we show that the rate of convergence depends on the "fractal" (i.e.~Minkowski-Bouligand) dimension of the boundary of the support of the graph limit. \newcor{Such a phenomenon was also reported in~\cite{medv} for a nonlocal nonlinear heat equation}. In Section~\ref{weightedgraphs}, we analyze networks on convergent weighted graph sequences. \newcor{Moreover, when the kernel and initial data belong to Lipschitz spaces, we also exhibit the convergence rate.} 

\paragraph{Notations}
For an integer $n \in \N^*$, we denote $[n]=\{1,\cdots,n\}$. For any set $\Omega$, $\cl{\Omega}$ is its closure, $\interop{\Omega}$ its interior and $\bd{\Omega}$ its boundary.



\section{Prerequisites on graphs}
\label{graphlimitss}
\subsection{Graph limits}

Let us start with reviewing some definitions and results from the theory of graph limits that we will need later since it is the key of our study of the discrete counterpart of the problem \eqref{neumann} on graphs. In our review, we follow considerably~\cite{graph,lovs}\newcor{, as presented in \cite{medv}}.

An undirected graph $G = \left ( V(G) , E(G) \right )$, where $V(G)$ stands for the set of nodes and $E(G) \subset V(G) \times V(G)$ denotes the edges set, without loops and parallel  edges is called simple.

\newcor{Let}  $G_{n} = \left ( V(G_{n}) , E(G_{n}) \right )$, $n \in \mathbb{N}$, be a sequence of dense, finite, and simple graphs, i.e; $\abs{E(G_{n})} = O (\abs{ V(G_{n})}^{2})$, where $\abs{.}$ denotes the cardinality of a set.

For two simple graphs $F$ and $G$, $\text{hom}(F , G)$ indicates the number of homomorphisms (adjacency-preserving maps) from $V (F )$ to $V (G)$. Then, it is worthwhile to normalize the homomorphism numbers and consider the homomorphism densities 
\[
t(F,G) = \frac{\text{hom} (F,G)}{\abs{V(G)}^{\abs{V(F)}}}.
\]
(Thus $t(F,G)$ is the probability that a random map of $V(F)$ into $V(G)$ is a homomorphism).
\begin{defi}\newcor{(cf.\cite{lovs})}
The sequence of graphs $\{G_{n}\}_n$ is called convergent if $ t(F,G_{n})$ is convergent for every simple graph $F$.
\label{convergence}
\end{defi}
 
\begin{rem}

Note that $t(F,G_{n} ) = O(1)$ if $\abs{E(G_{n})} = O (\abs{ V(G_{n})}^{2})$  so that this definition is meaningful only for sequences of dense graphs. In the theory of graph limits, convergence in Definition~\ref{convergence} is called left-convergence. Since this is the only convergence of graph sequences that we use, we would refer to the left-convergent sequence as convergent (see \cite[Section~2.5]{borg}).
\end{rem}

Convergent graph sequences have a limit object, which can be represented as a measurable \newcor{symmetric} function $K: \O^{2} \to \O$, here $\O$ stands for $[0,1]$. Such functions are called graphons.

Let $\mathcal{K}$ denote the space of all bounded measurable functions  $K : \O^2 \to  \R$ such that $K(x,y)=K(y,x) $ for all $x,y \in [0,1]$. We also define $\mathcal{K}_{0} =\{ K \in \mathcal{K} : 0 \le K \le 1\}$ the set of all graphons.
\begin{prop}[{\cite[Theorem~2.1]{graph}}]
 For every convergent sequence of simple graphs, there is $K \in  \mathcal{K}_0$ such that 
 \begin{equation}
 t(F,G_{n}) \to t(F,K) := \int_{\O^{\abs{V(F)}}}  \prod_{(i,j) \in E(F)} K(x_{i},x_{j}) dx.
 \label{lim}
 \end{equation}
for every simple graph $F$. Moreover, for every $K \in \mathcal{K}_{0} $, there is a sequence of graphs $\{G_{n}\}_n$ satisfying~\eqref{lim}.
\end{prop}
 
Graphon $K$ in \eqref{lim} is the limit of the convergent sequence $\{G_{n}\}_n$. It is uniquely determined up to measure-preserving transformations in the following sense: for every other limit function $K' \in \mathcal{K}_{0} $, there are measure-preserving map $\phi , \psi: \O \to \O$ such that $K(\phi(x), \phi(y)) = K' (\psi(x), \psi(y))$ (see \cite[Theorem~2.1]{graph}).
 
Indeed, every finite simple graph $G_{n}$ such that $V(G_{n}) = [n]$ can be represented by a function $K_{G_{n}} \in \mathcal{K}_{0}$ 
\begin{equation*}
 K_{G_{n}} (x,y) = 
 \begin{cases}
 1 \quad \text{if}\quad (i,j) \in E(G_{n}) \quad \text{and} \quad (x,y) \in [\frac{i-1}{n} , \frac{i}{n}[ \times [\frac{j-1}{n} , \frac{j}{n}[, \\
0 \quad \text{otherwise}.
 \end{cases}
\end{equation*}
Hence, geometrically, the graphon $K$ can be interpreted as the limit of $K_{G_{n}}$ for the standard (called the cut-norm) 
\[
 \norm{K}_{ \square} := \sup_{S,T \in \mathcal{L}_{\O}} \abs{\int_{S \times T} K(x,y) dx dy},
\]
 where $K  \in L^{1}( \O^2) $ and $\mathcal{L}_{\O}$ stands for the set of all Lebesgue measurable subsets of $\O$.
 Since for any $K  \in L^{1}( \O^2) $ 
\[
  \norm{K}_{ \square}   \le  \norm{K}_{ L^1(\O^2)},
\]
  convergence of $\{ K_{G_{n}}\}$ in the $L^1$-norm implies the convergence of the graph sequence $\{ G_{n}\}_n$ (\cite[Theorem~2.3]{graph}).

We finish this section by giving an example of convergent graph sequences that is very useful in practice.
\begin{exe} (see \cite{lovs}) The Erd\"os-Renyi graphs. Let $p \in ]0,1[$ and consider the sequence of random graphs $G(n,p) = \left ( V(G(n,p)) , E(G(n,p)) \right )$ such that $V(G(n,p)) = [n]$ and the probability \linebreak$\Pr \{ (i,j) \in E(G(n,p))\}=p $ for any $(i,j) \in [n]^2$. Then for any simple graph $F$, $t(F, G(n,p))$ is convergent with probability $1$ to $p^{\abs{E(F)}}$ as $n \to \infty$ \cite{borg}.
\end{exe}  

\subsection{Types of graph sequences}
 \label{subsec:types}

\newcor{The graph models presented below were constructed in~\cite{medv}}.
\subsubsection{Simple graph sequences}
We fix $n \in \N$, divide $ \O$ into $n$ intervals 
\[
\O_1^{(n)} =  \left [0, \frac{1}{n} \right [ , \O_2^{(n)} =  \left [\frac{1}{n}, \frac{2}{n} \right [,\ldots,\O_j^{(n)} =  \left [\frac{j-1}{n}, \frac{j}{n} \right [,\ldots,\O_n^{(n)} =  \left [ \frac{n-1}{n} , 1\right [,
\]
and let  $\mathcal{Q}_n$ denote the partition of $\O$, $\mathcal{Q}_n = \{ \O_i^{(n)}, i \in [n]\}$. Denote $\O_{ij}^{(n)} := \O_{i}^{(n)} \times \O_{j}^{(n)}$.

 We consider first the case of a sequence of simple graphs converging to $\{ 0,1\}$ graphon.

Briefly speaking, we define a sequence of simple graphs $G_n = \left (V(G_n),E(G_n) \right )$ such that $V(G_n) = [n]$ and 
\[
E(G_n) = \ens{(i,j) \in [n]^2}{ \O_{ij}^{(n)} \cap \cl{\supp(K)} \ne \emptyset},
\]
where 
\begin{equation}
\supp(K) = \ens{(x,y) \in \O^2}{K(x,y) \ne 0}.
\label{support}
\end{equation}
As we have mentioned before, the kernel $K$ represents the corresponding graph limit, that is the limit as $n \to \infty$ of the function $K_{G_n} : \O^2 \to \{0,1\}$ such that 
\[
 K_{G_n}(x,y)   = \left\{\begin{array}{lr}
        1, \quad \text{if} \quad (i,j) \in E(G_n) \quad \text{and} \quad (x,y) \in \O_{ij}^{(n)},
        \\
      0  \quad \text{otherwise}.
        \end{array}
        \right.
\] 
As $n \to \infty $, $\{K_{G_n}\}_n$ converges to the $\{0,1\}$-valued mapping $K (\cdot, \cdot)$ whose support is defined by~\eqref{support}.

\subsubsection{Weighted graph sequences} 
\newcor{We now review a more general class of graph sequences. We consider two sequences of weighted graphs generated by a given graphon $K$.}

Let $K: \O^2 \to [a,b]$ $a,b > 0$, be a symmetric measurable function which will be used to assign weights to the edges of the graphs considered bellow, we allow only positive weights.

Next, we define the quotient of $K $ and $\mathcal{Q}_n$ denoted $K / \mathcal{Q}_{n} $ as a weighted graph with $n$ nodes  
\[
K / \mathcal{Q}_{n} = \left ( [n] , [n] \times [n] , \hat{K}_{n}  \right ).
\]

As before, weights $(\k)_{ij}$ are obtained by averaging $K$ over the sets in $\mathcal{Q}_n$ 
\begin{equation}
\label{eq:kweight}
(\k)_{ij} = n^2 \int_{\O^{(n)}_{i} \times \O^{(n)}_{j}}  K(x,y) dx dy.
\end{equation}


The second sequence of weighted graphs is constructed as follows 
\[
\mathbb{G} (X_{n},K) = \left ([n] , [n] \times [n] , \breve{K_{n}} \right ),
\]
where 
\begin{equation}
X_{n} = \left \{  \frac{1}{n} ,  \frac{2}{n}, \cdots ,  \frac{n}{n}  \right \}, \hspace{0.5cm}  (\breve{K_{n}})_{ij} = K \left ( \frac{i}{n} ,  \frac{j}{n} \right ).
\label{collocation}
\end{equation}

\section{Existence and uniqueness of a solution}
\label{existenceuniqueness}
The main result of existence and uniqueness of a global solution, that is, a solution on $[0, T]$ for $T>0$ is stated in the following theorem.

\begin{theo}
\label{theo:exist}
Suppose $p \in ]1,+\infty[$ and let $g \in L^p(\O)$.
\begin{enumerate}[label=(\roman*)]
\item For any $T > 0$, there exists a unique strong solution in $[0,T]$ of \eqref{neumann}.
\item Moreover, for $q \in [1,+\infty]$, if $g_i \in L^q(\O)$, $i=1,2$, and $u_i$ is the solution of \eqref{neumann} with  initial condition $g_i$, then
\begin{equation}
\norm{u_1(t) - u_2(t)}_{L^q(\O)} \le \norm{g_1 - g_2}_{L^q(\O)}, \quad \forall t \in [0,T] .
\label{contraction}
\end{equation}

\end{enumerate}
\label{existence}
\end{theo}
\begin{rem}
For $p \in [1,+ \infty]$, taking the initial data in $L^p(\O) $, one can show existence and uniqueness of a mild but not a strong solution as  $L^1(\O)$ and $L^{\infty}(\O)$ are not reflexive spaces and thus do not have the Radon-Nikodym property (see  \cite[Theorem~A.29 and Proposition~A.35]{vaillo}).
\end{rem}
 The proof of Theorem~\ref{theo:exist} is an extension of that of \cite[Theorem~6.8]{vaillo} to the case of a symmetric, nonnegative and bounded kernel $K$ as in our setting (see \cite[Remark~6.9]{vaillo}). For this, we only need to show the corresponding versions of \cite[Lemmas~6.5 and 6.6]{vaillo} (which are stated there without a proof). See Section~\ref{subsec:prooftheoexist} for details.

\section{Consistence of the time-continuous problem} 
\label{timecontinuous}
We begin our study by giving a general consistency result from which we shall extract particular consistency bounds for every specific model of convergent graph sequences that we have introduced in section~\ref{subsec:types}. To do this, let us consider the following Cauchy problem with Neumann boundary conditions as \eqref{neumann}
\begin{equation}\tag{\textrm{$\mathcal{P}_n$}}
\begin{cases}
 \frac{\partial }{\partial t} u_{n}(x,t) = -\An(u_{n}(x,t)), & (x,t) \in \O \times ]0,T]
  \\
  u_{n}(x,0)= g_{n}(x), \quad x \in \O .
\end{cases}
\label{general}
\end{equation}
Though not needed in this section, the use of the subscript $n$ is a matter of notation and emphasizes the fact that $K_n$ and $g_n$ depend on the parameter $n$. This will be clear in the application to graphs (Section~\ref{applications}).  

Now we state and prove our main uniform convergence theorem.
\begin{theo}
Suppose $p \in ]1, +\infty[$, $g, g_{n} \in L^{\infty}(\O)$ and $\K , \Kn$ are measurable, symmetric and bounded mappings. Then \eqref{neumann} and \eqref{general} have unique solutions, respectively, $u$ and $u_n$. Moreover the following \newcor{hold}. 
\begin{enumerate}[label=(\roman*)]
\item We have the error estimate
\begin{equation}
\norm{u- u_n}_{C(0,T;L^p(\O))} \le C \left (  \norm{g - g_n}_{L^p(\O)} + \norm{\K - \Kn}_{L^p(\O^2)}  \right ),
\label{stability}
\end{equation}
where the constant $C$ is independent of $n$. 
\item
\label{scheffe}
 Moreover, if $g_n \to g$ and $\Kn \to \K$ as $n \to \infty$, almost everywhere on $\O$ and $\O^2$, respectively, then 
\[
\norm{u- u_n}_{C(0,T;L^p(\O))} \underset{{n \to \infty}}{\longrightarrow} 0.
\]
\end{enumerate}
\label{main}
\end{theo}
\bpf{}
In the proof, $C_i$ is any absolute constant independent of $n$ (but may depend on $p$). Existence and uniqueness of the solutions $u$ and $u_n$ in the sense of Definition~\ref{def:sol} is a consequence of Theorem~\ref{existence}.
\begin{enumerate}[label=(\roman*)]
\item For $1<p<+ \infty$, we define the function 
  \begin{equation*}
  \begin{split}
  \F :   x \in \mathbb{R} \mapsto \abs{x}^{p-2} x = \sign(x) {\abs{x}}^{p-1}.
  \end{split}
 \end{equation*}
Denote $\xi_n(x,t) = u_n(x,t) - u(x,t)$, by subtracting \eqref{neumann} from \eqref{general}, we have a.e.
 \begin{equation}
 \begin{split}
 \frac{\partial \xi_n(x,t)}{\partial t }  &= \int_{\O} \Kn(x,y) \{ \F(u_n(y,t) - u_n(x,t)) - \F(u(y,t)-u(x,t))\} dy\\
 &+ \int_{\O} (\Kn(x,y) - \K(x,y)) \F(u(y,t)-u(x,t)) dy.
 \end{split}
 \label{deriv}
 \end{equation}
  Next, we multiply both sides of \eqref{deriv} by $\F(\xi_n(x,t))$ and integrate over $\O$ to get
  \begin{equation}
  \begin{split}
  \frac{1}{p} \int_{\O} \frac{\partial}{\partial t} \abs{\xi_n(x,t)}^p dx 
  &=  \int_{\O^2} \Kn(x,y) \{ \F(u_n(y,t) - u_n(x,t)) - \F(u(y,t)-u(x,t))\} \F(\xi_n (x,t)) dx dy \\
   &+ \int_{\O^2} ( \Kn(x,y) - \K(x,y) ) \F(u(y,t)-u(x,t)) \F(\xi_n (x,t)) dx dy.
  \end{split}
  \label{derivv}
  \end{equation}
  We estimate the first term on the right-hand side of \eqref{derivv} using the fact that $\Kn$ is bounded so that there exists a positive constant $M$ independent of $n$, such that,  $\norm{\Kn}_{L^{\infty}(\O^2)} \le M$,
\begin{align*}
& \abs{\int_{\O^2} \Kn(x,y) \{ \F(u_n(y,t) - u_n(x,t)) - \F(u(y,t)-u(x,t))\} \Psi(\xi_n (x,t)) dx dy } \\
& \le M \int_{\O^2} \abs { \F(u_n(y,t) - u_n(x,t)) - \F(u(y,t)-u(x,t)) } \abs{\xi_n(x,t) }^{p-1} dx dy .
\end{align*}
Now, applying Lemma~\ref{cor:genmvtlp} with $a=u_n(y,t) - u_n(x,t)$ and $b=u(y,t)-u(x,t)$ (without loss of generality we assume that $b > a$), we get
\begin{align}
& \int_{\O^2} \abs { \F(u_n(y,t) - u_n(x,t)) - \F(u(y,t)-u(x,t)) } \abs{\xi_n (x,t) }^{p-1} dx dy \nonumber\\
\label{eq:bndmvtlp}
& \le (p-1) \int_{\O^2}  \abs{\xi_n(y,t) - \xi_n(x,t)} \abs{\eta(x,y,t)}^{p-2}  \abs{\xi_n (x,t) }^{p-1} dx dy,
\end{align}
where $\eta(x,y,t) \neq 0$ is an intermediate value between $a$ and $b$. \newcor{As we have supposed that $g \in L^{\infty} (\O)$ and $g_n \in L^{\infty}(\O)$, and as $|\O|$ is finite, so that $L^{\infty} (\O) \subset L^p(\O) $, we deduce from \eqref{contraction} in Theorem~\ref{existence} that for any $(x,y) \in \O^2$ and $t \in [0,T]$, we have for $p \geq 2$} 
\begin{equation*}
\begin{aligned}
\abs{\eta(x,y,t)}^{p-2} \le \abs{u(y,t)-u(x,t)}^{p-2} &\le \pa{2\norm{u(t)}_{L^{\infty}(\O)}}^{p-2} \leq \pa{2\norm{g}_{L^{\infty}(\O)}}^{p-2} < +\infty .
\end{aligned}
\end{equation*} 
For $p \in ]1,2[$, since $\inf_{(x,y,t) \in \O^2 \times [0,T]} |\eta(x,y,t)| =  C' > 0$, we have
\begin{equation*}
\begin{aligned}
\abs{\eta(x,y,t)}^{p-2} \le C'^{p-2} < +\infty .
\end{aligned}
\end{equation*}
Thus, letting $C_1 = \pa{\max\pa{2\norm{g}_{L^{\infty}(\O)}^{p-2},C'^{p-2}}}$, we have
\begin{equation}
\begin{aligned}
 \abs{\eta(x,y,t)}^{p-2} \le C_1 .
\end{aligned}
\label{major}
\end{equation}
Inserting \eqref{major} into \eqref{eq:bndmvtlp}, and then using the H\"older and triangle inequalities, it follows that  
\begin{equation}
  \begin{split}
 &M \int_{\O^2} \abs {\F(u_n(y,t) - u_n(x,t)) - \F(u(y,t)-u(x,t)) } \abs{\xi_n (x,t) }^{p-1} dx dy \\
 &\le M (p-1)C_1 \int_{\O^2} \abs{\xi_n(y,t) - \xi_n(x,t)}  \abs{\xi_n (x,t)}^{p-1} dx dy \\
 &= C_2   \int_{\O^2} \abs{\xi_n(y,t) - \xi_n(x,t)}  \abs{\xi_n (x,t)}^{p-1} dx dy \\
 & \le C_2   \left ( \int_{\O^2} \abs{\xi_n(y,t) - \xi_n(x,t)}^{p} dx dy \right )^{\frac{1}{p}} \times \left ( \int_{\O} \abs{\xi_n(x,t)}^{p} dx \right )^{\frac{p-1}{p}}   \\
 &\le 2 C_2 \norm{\xi_n(t)}_{L^p(\O)}^p.
  \end{split}
  \label{est}
\end{equation}
We bound the second term on the right-hand side of \eqref{derivv} as follows
\begin{equation}
\begin{split}
& \abs { \int_{\O^2} (\Kn(x,y) - \K(x,y)) \F(u(y,t) - u(x,t) ) \F(\xi_n (x,t)) dx dy }\\
&=\abs {  \int_{\O^2} (\Kn(x,y) - \K(x,y) ) \times \sign (u(y,t) - u(x,t) ) \abs{u(y,t) - u(x,t) }^{p-1}  \F(\xi_n (x,t)) dx dy } \\
&\le 2^{p-1} \norm{u(t)}_{L^{\infty}(\O)}^{p-1} \abs {  \int_{\O^2} \abs{\Kn(x,y) - \K(x,y)}  \abs{\xi_n(x,t)}^{p-1}   dx dy } \\
&\le 2^{p-1} \norm{u(t)}_{L^{\infty}(\O)}^{p-1}  \left ( \int_{\O} \abs{\xi_n(x,t)}^{p} dx \right )^{\frac{p-1}{p}} \times   \left ( \int_{\O^2} \abs{\Kn(x,y) - \K(x,y)}^p dx dy  \right )^{\frac{1}{p}} \\
&\le 2 C_3  \norm{\xi_n(t)}_{L^{p}(\O)}^{p-1} \norm{\Kn - \K}_{L^{p}(\O^2)}.
\end{split}
\label{esti}
\end{equation}
Bringing together \eqref{est} and \eqref{esti}, and using standard arguments to switch the derivation and integration signs (Leibniz rule), we have
\begin{equation}
\frac{d}{dt} \norm{\xi_n(t)}_{L^p(\O)}^p \le 2pC_2 \norm{\xi_n(t)}_{L^p(\O)}^p + 2pC_3 \norm{\Kn - \K}_{L^{p}(\O^2)} \norm{\xi_n(t)}_{L^p(\O)}^{p-1}.
\label{resssi}
\end{equation}
Let $\varepsilon > 0$ be arbitrary but fixed, and set 
\[
\psi_{\varepsilon} (t) = \pa{\norm{\xi_n(t)}_{L^p(\O)}^p + \varepsilon}^{1/p}.
\]
By \eqref{resssi},
\begin{equation}
\frac{d}{dt}  \psi_{\varepsilon} (t)^p \le 2pC_2 \psi_{\varepsilon} (t)^p + 2pC_3  \norm{\Kn - \K}_{L^{p}(\O)} \psi_{\varepsilon} (t)^{p-1}.
\label{ress}
\end{equation}
 Since $\psi_{\varepsilon} (t)$ is positive on $[0,T]$, from \eqref{ress}, we have 
\[
\frac{d}{dt}  \psi_{\varepsilon} (t) \le 2 C_2  \psi_{\varepsilon} (t) + 2 C_3 \norm{\Kn - \K}_{L^{p}(\O^2)} ,  \hspace{0.5cm} t \in [0,T].
\]
We apply \newcor{ Gronwall's} inequality for $\psi_{\varepsilon} (t)$ on $[0,T]$ to get 
\begin{equation}
\sup_{t \in [0,T]} \psi_{\varepsilon} (t) \le \left (  \psi_{\varepsilon} (0) + 2C_3 T \norm{\Kn - \K}_{L^{p}(\O^2)}\right ) \exp \{2 C_2 T\}.
\label{gron}
\end{equation}
Since $\varepsilon > 0$ is arbitrary, \eqref{gron} implies 
\begin{equation}
\sup_{t \in [0,T]} \norm{\xi_n(t)}_{L^p(\O)} \le \left (  \norm{g-g_n}_{L^p(\O)} + 2C_3 T  \norm{\Kn - \K}_{L^{p}(\O^2)} \right ) \exp \{2 C_2 T\}.
\label{resss}
\end{equation}
The desired result holds.
\item Since $g_n, g \in L^\infty(\O) \subset L^p(\O)$ and $|\O|$ is finite, the dominated convergence theorem implies that $\lim_{n \to +\infty} \norm{g_n}_{L^p(\O)} = \norm{g}_{L^p(\O)}$. The same reasoning applies to $\Kn$ and $\K$. Passing to the limit in \eqref{stability} and using the Scheff\'e-Riesz theorem (see \cite[Lemma~2]{riesz}), we get the claim.
\end{enumerate}
\epf
\begin{rem}
Observe that, \newcor{since $|\O|$ is finite}, we have the classical inclusion $L^p(\O) \subset L^2 (\O)$ for $p \ge 2$, which leads to the following bound 
\[
\norm{u- u_n}_{C(0,T; L^2(\O))} \le \abs{\O}^{\frac{1}{2} - \frac{1}{p}} \norm{u- u_n}_{C(0,T; L^p(\O))} = \norm{u- u_n}_{C(0,T; L^p(\O))},
\]
as $\abs{\O}=1$. For $p \in ]1,2]$, we have, thanks to Lemma~\ref{inclusion}, boundedness of the solutions and Jensen inequality,
\[
\norm{u- u_n}^2_{C(0,T; L^2(\O))} = O\pa{\norm{u- u_n}^p_{C(0,T; L^p(\O))}} = O\pa{\norm{g - g_n}^p_{L^p(\O)} +\norm{\K - \Kn}^p_{L^p(\O^2)}} .
\]
In summary, there is also convergence with respect to the $L^2$-norm.
\end{rem}

\section{Consistence of the time-discrete problem}
\label{timediscrete}

\subsection{Forward Euler discretization}\label{subsec:Pdf}
We now consider the following time-discrete approximation of~\eqref{neumann}, the forward Euler discretization applied to~\eqref{general}. \newcor{For that, let us consider a partition (not necessarily uniform) $\{t_h\}_{h=1}^{N}$ of the time interval $[0,T]$. Let $\tauhm := \abs{  t_h -  t_{h-1}}$ and the maximal size $\tau = \max\limits_{h \in [N]} \tauh $, and denote $\h (x) := u_n(x,t_h)$}. Then, consider

\begin{equation}\tag{\textrm{$\mathcal{P}^{f}_{n,\tau}$}}
\begin{cases}
\displaystyle{\frac{\h(x) - \H(x)}{\tauhm}} = -\An(\H(x)), & x \in \O, h \in [N],\\
u_n^{0}(x)= \g(x), \quad x \in \O .
\end{cases}
\label{discretisation}
\end{equation}
Before turning to the consistency result, one may wonder whether \eqref{discretisation} is well-posed. In the following result, we show that for $p \in ]1,+\infty[$, and starting from $\g \in L^\infty(\O)$, there exists a unique weak accumulation point to the iterates of \eqref{discretisation}. In turn, in the case of practical interest where the problem is finite-dimensional (in fact Euclidean case) as for the application to graphs (see Section~\ref{applications}), we do have existence and uniqueness. Recall the function $F_p$ from \eqref{eq:Fp}.
\begin{lem} 
Consider problem~\eqref{discretisation}. Assume that $\g \in L^\infty(\O)$. Let $\tauh =\displaystyle{ \frac{\alpha_h}{\max(\norm{\An (\h)}_{L^2(\O)},1)}}$, and suppose that $\sum\limits_{h=1}^{+\infty} \alpha_{h} = + \infty$ and $\sum\limits_{h=1}^{+\infty} \alpha_{h}^2 < + \infty$. Then, the iterates of problem~\eqref{discretisation}, starting from $\g$, have a unique weak accumulation point $u^{\star}$. Moreover, there are constants $\beta,\varepsilon > 0$ such that
\[
\min_{0 \leq i \leq h} F_p(u_n^{i}) - F_p(u^\star) \leq \max(\beta,1)\frac{\varepsilon^2+\sum_{i=0}^h \alpha_i^{2}}{2\sum_{i=0}^h \alpha_i} .
\]
\label{lem:existPdf}
\end{lem}
\begin{rem}
\begin{enumerate}[label=(\alph*)]
\item Our condition on the time-step $\tauh$ is reminiscent of the subgradient method. It can be seen as a non-linear CFL-type condition which depends on the data since $\An$ is not Lipschitz-continuous but only locally so, hence the dependence of $\tauh$ on $\norm{\An (\h)}_{L^2(\O)}$.
\item The rate of convergence on $F_p$ depends on the choice of $\acc{\alpha_h}_h$. If one performs $N$ steps on the interval $[0,T]$, one can take
\[
\alpha_h = \frac{\varepsilon}{(N+1)^{1/2+\nu}} , h = 0, \ldots, N, \qwithq \nu \in ]0,1/2[ ,
\]
which entails a convergence rate of $\tfrac{\max(\beta,1)\varepsilon^2}{(N+1)^{1/2-\nu}}$. The smaller $\nu$ the faster the rate.
\end{enumerate}
\end{rem}
Before proving Lemma~\ref{lem:existPdf} recall the definition of the subdifferential. Let $F: L^2(\O) \to \R \cup \acc{+\infty}$ be a proper lower-semicontinuous and convex function. The subdifferential of $F$ at $u \in L^2(\O)$ is the set-valued operator $\partial F: L^2(\O) \to 2^{L^2(\O)}$ given by
\[
\partial F(u) = \ens{\eta \in L^2(\O)}{F(v) -F(u) \ge \pds{\eta}{u-v}, \quad \forall v \in L^2(\O)},
\]
\newcor{where $\pds{.}{.}$ denotes the inner product in $L^2(\O)$}. 

Moreover, $F$ is G\^ateaux differentiable at $u \in \dom(F)$ if and only if $\partial F(u)$ is a singleton with the gradient $\nabla F(u)$ as its unique element~\cite[Corollary~17.26]{bauschke2011convex}.

\newcor{\bpf{}} 
Since $p > 1$, we consider in the Hilbert space $L^2(\O)$ the subdifferential $\partial F_p$ whose graph is in $L^2(\O) \times L^2(\O)$. It is immediately seen that $F_p$ is convex and G\^ateaux-differentiable, and thus $\partial F_p(u)=\acc{\An(u)}$. Moreover, it is maximal monotone (or equivalently $m$-accretive on $L^2(\O)$), see \cite[p.~198]{vaillo}. Consequently, using that $\g \in L^{\infty}(\O) \subset L^2(\O)$, and so is $\h$ by induction, a solution to \eqref{discretisation} coincides with that of 
\begin{equation*}
\begin{cases}
\h(x) \in \H(x) - \tauhm \eta^{h-1}, & \eta^{h-1} \in \partial F_p(\h) \\
u_n^{0}(x)= \g(x), \hspace{1cm} x \in \O ,
\end{cases}
\end{equation*}
i.e. the subgradient method with initial point $\g$. Observe that $(\partial F_p)^{-1} (0) \ne \emptyset$ ($0$ is in it). Thus with the prescribed choice of $\tauh$, we deduce from \cite[Theorem~1]{Alber98} that the sequence of iterates $\h$ has a unique weak accumulation $u^\star \in (\partial F_p)^{-1} (0)$.

The claim on the rate is classical~\footnote{See e.g. \cite[Theorem~3.2.2]{nesterov2004introductory} in finite dimension with a slightly different normalization of the step size $\tauh$.}. We here provide a simple and self-contained proof. Since $F_p$ is continuous and convex on $L^2(\O)$, it is locally Lipschitz continuous~\cite[Theorem~8.29]{bauschke2011convex}. Moreover, the sequence $\acc{\h}_h$ is bounded, and hence, $\exists \varepsilon > 0$ such that $\norm{\h-u^\star}_{L^2(\O)} \leq \varepsilon$, $\forall h \geq 0$. In turn, $F_p$ is Lipschitz continuous around $u^\star$ with Lipschitz constant, say $\beta$. Denote $\rh=\h-u^{\star}$. We have 
\begin{equation*}
\begin{aligned}
\norm{\rh}_{L^2(\O)}^2 
&= \norm{\rH  - \tauhm \eta^{h-1}}_{L^2(\O)}^2\\
& = \norm{\rH}_{L^2(\O)}^2 - 2\tfrac{\alhm}{\max\pa{\norm{\eta^{h-1}}_{L^2(\O)},1}}\pds{\eta^{h-1}}{\rH} + \alhm^2 \\
&\leq \norm{\rH}_{L^2(\O)}^2 - 2\tfrac{\alhm}{\max\pa{\norm{\eta^{h-1}}_{L^2(\O)},1}}\pa{F(\H))-F_p(u^\star)} + \alhm^2,
\end{aligned}
\end{equation*}
where we used the subdifferential inequality above to get that
\[
F_p(u^\star) \geq F(\H) - \pds{\eta^{h-1}}{\rH} .
\]
Summing up these inequalities we obtain
\begin{align*}
2\sum_{i=0}^h \alpha_i \pa{F_p(u_n^{i}) - F_p(u^\star)} \leq \max(\beta,1) \pa{\norm{r_n^{0}}_{L^2(\O)}^2 + \sum_{i=0}^h \alpha_i^2} ,
\end{align*}
whence we deduce
\[
\min_{0 \leq i \leq h} F_p(u_n^{i}) - F_p(u^\star) \leq \max(\beta,1)\frac{\varepsilon^2+\sum_{i=0}^h \alpha_i^{2}}{2\sum_{i=0}^h \alpha_i} .
\]
\epf

Since the aim is to compare the solutions of problems~\eqref{neumann} and~\eqref{discretisation}, the solution of~\eqref{discretisation} being discrete, so that it is convenient to introduce an intermediate model which is the continuous extension of the discrete problem using the discrete function $u_n(x) = (u_n^1(x), \cdots,  u_n^N(x))$. Therefore, we consider a time-continuous extension of $\h$ obtained by a time linear interpolation as follows 
\begin{equation}
\U(x,t) = \frac{t_{h} - t}{\tauhm}  \H(x)+ \frac{t-t_{h-1}}{\tauhm} \h(x) ,\quad t \in ]t_{h-1},t_h] , \quad x \in \O,
\label{convex}
\end{equation}
and a time piecewise constant approximation
\begin{equation}
\M (x,t) = \sum_{h=1}^N \H(x) \chi_{]t_{h-1},t_h]} (t).
\label{piecewise}
\end{equation}
Then, by construction of $\U(x,t)$ and $\M(x,t)$, we have the following evolution problem 
\begin{equation}
\begin{cases}
\frac{\partial}{\partial t} \U(x,t)=-\An(\M(x,t)), & (x,t) \in \O \times ]0,T]\\
\U(x,0) = g_n^0(x), \quad x \in \O .
\end{cases}
\label{relation}
\end{equation}

\begin{lem} 
Assume that $\g \in L^\infty(\O)$. Let $\U$ and $\M$ be the functions defined in~\eqref{convex} and~\eqref{piecewise}, respectively, then 
\begin{equation}
\norm{\M (t)- \U(t)}_{L^p(\O)} = O(\tau), \quad \quad t\in [0,T].
\end{equation}
\label{intermediaire}
\end{lem}
\bpf{}
It is easy to see that for  $t \in ]t_{h-1},t_h]$,
\begin{align*}
\norm{\M(t) - \U(t)}_{L^p(\O)} \leq (t_h - t) \norm{\frac{\h -\H}{\tauhm}}_{L^p(\O)} 
&\leq \tau \norm{\frac{\h-\H}{\tauhm}}_{L^p(\O)} = \tau \norm{\An(\H)}_{L^p(\O)} \\
&\leq \tau \norm{\An(\H)}_{L^\infty(\O)} \leq \tau 2^{p-1}\norm{\H}_{L^\infty(\O)}^{p-1} .
\end{align*}
By induction, for all $h \geq 1$, we have (see Lemma~\ref{lem:existPdf})
\[
\norm{\h}_{L^\infty(\O)} \leq \norm{\H}_{L^\infty(\O)} + \alpha 2^{p-1}\norm{\H}_{L^\infty(\O)}^{p-1} < +\infty , 
\]
where $\alpha  = \sup\limits_{h \geq 1} \alpha_h < +\infty$. Since $t$ is arbitrary, we obtain a global estimate for all $t \in [0,T]$. 
\epf

\begin{theo}
Suppose $p \in ]1, +\infty[$, $g, \g\in L^{\infty}(\O)$ and $\K, \Kn$ are measurable, symmetric and bounded mappings.

Let $u$ be the unique solution of problem~\eqref{neumann}, and $\U$ is built as in~\eqref{convex} from the time-discrete approximation $\H$ defined in~\eqref{discretisation}. Then 
\begin{equation}
\norm{u- \U}_{C(0,T;L^p(\O))} \le C \left (  \norm{g_n - \g}_{L^p(\O)} +\norm{g - g_n}_{L^p(\O)} + \norm{\K - \Kn}_{L^p(\O^2)}  \right ) +O(\tau),
\label{stability1}
\end{equation}
where the constant $C$ is independent of $n$. 
\label{totally}
\end{theo}
\bpf{}
We follow the same lines as in the proof of Theorem~\ref{main}. Denote $\q(x,t) = \U(x,t) - u_n(x,t)$ and $\s(x,t) = \M(x,t) - u_n(x,t)$. We thus have a.e.
\begin{equation}
 \frac{\partial \q}{\partial t }  = \int_{\O} K_{n}(x,y) \{ \F(\M(y,t) - \M(x,t)) - \F(u_n(y,t)-u_n(x,t))\} dy.
  \label{derivvv}
\end{equation}
Next, we multiply both sides of \eqref{derivvv} by $\F(\q (x,t))$ and integrate over $\O$ using the relation~\eqref{relation} to get
  \begin{equation}
  \frac{1}{p} \int_{\O} \frac{\partial}{\partial t} \abs{\q(x,t)}^p dx 
  =  \int_{\O^2} \Kn(x,y) \{ \F(\M(y,t) - \M(x,t)) - \F(u_n(y,t)-u_n(x,t))\} \F(\q) (x,t) dx dy.
   \label{derivvv1}
  \end{equation}
  
Similarly to the proof of Theorem~\ref{main}, we bound the term on the right-hand side of \eqref{derivvv1} using the fact that $K_n$ is bounded, then applying Corollary~\ref{cor:genmvtlp} between $ \M(y,t) - \M(x,t)$ and $ u_n(y,t)-u_n(x,t)$, inequality \eqref{major}, and finally using H\"older and triangle inequalities. Altogether, this yields 
\begin{equation}
\begin{split}
&\abs{\int_{\O^2} \Kn(x,y) \{ \F(\M(y,t) - \M(x,t)) - \F(u_n(y,t)-u_n(x,t))\} \F(\q) (x,t) dx dy } \\
&\le  C_2   \int_{\O^2} \abs{\s(y,t) - \s(x,t)}  \abs{\q (x,t)}^{p-1} dx dy \\
& \le C_2   \left ( \int_{\O^2} \abs{\s(y,t) - \s(x,t)}^{p} dx dy \right )^{\frac{1}{p}} \times \left ( \int_{\O} \abs{\xi_n(x,t)}^{p} dx \right )^{\frac{p-1}{p}}   \\
 &\le 2 C_2 \norm{\s(t)}_{L^p(\O)} \norm{\q(t)}_{L^p(\O)}^{p-1}.
  \end{split}
  \label{holder}
  \end{equation}
  
By virtue of Lemma~\ref{intermediaire} and the triangle inequality for $\s(\cdot, \cdot)$, there exists a positive constant $C'$ such that  
  \begin{equation}
  \begin{aligned}
  \norm{\M(t) - u_n(t)}_{L^p(\O)} &\le \norm{\M (t)- \U(t)}_{L^p(\O)} + \norm{\U(t)- u_n(t)}_{L^p(\O)}\\
 &\le C' \tau + \norm{\q(t)}_{L^p(\O)}.
  \end{aligned}
  \label{inter}
  \end{equation}
  Hence, bringing together \eqref{holder} and \eqref{inter}, we obtain
  \begin{equation}
\frac{d}{dt} \norm{\q(t)}_{L^p(\O)}^p \le 2pC_2 \norm{\q(t)}_{L^p(\O)}^p + 2 p C' \tau \norm{\q(t)}_{L^p(\O)}^{p-1}.
\end{equation}

Arrived at this stage, we proceed in the same way using the Gronwall's lemma as in the proof of Theorem~\ref{main}, to get 
\begin{equation}
\sup_{t \in [0,T]} \norm{\q(t)}_{L^p(\O)} \le \left (  \norm{\g-g_n}_{L^p(\O)} + 2 C' T \tau \right ) \exp \{2 C_2 T\}.
\end{equation}
Then, 
\begin{equation}
\norm{\U- u_n}_{C(0,T;L^p(\O))} \le C  \norm{\g - g_n}_{L^p(\O)} + C^{''} \tau.
\label{approx}
\end{equation}
Using the triangle inequality and \eqref{stability} in Theorem~\ref{main}, we get 
\begin{equation}
\begin{aligned}
\norm{\U - u}_{C(0,T;L^p(\O))} &\le \norm{\U- u_n}_{C(0,T;L^p(\O))}  + \norm{u_n - u }_{C(0,T;L^p(\O))}\\
& \le C^{''} \tau + C \left ( \norm{\g - g_n}_{L^p(\O)}  + \norm{g - g_n}_{L^p(\O)} + \norm{\K - \Kn}_{L^p(\O^2)}  \right ).
\end{aligned}
\end{equation}
\epf

\subsection{Backward Euler discretization}\label{subsec:Pdb}
Our result in Theorem~\ref{totally} also holds when we deal with the backward Euler discretization
\begin{equation}\tag{\textrm{$\mathcal{P}_{n,\tau}^b$}}
\begin{cases}
\displaystyle{\frac{\h(x) - \H(x)}{\tauhm}} = -\An(\h(x)), & x \in \O, h \in [N], \\
u^{0}(x)= \g(x), \quad x \in \O,
\end{cases}
\label{discretisationb}
\end{equation}
which can also be rewritten as the implicit update
\begin{equation*}
\begin{cases}
\h(x) = J_{\tauhm\An}(\H)(x), & x \in \O, h \in [N],\\
u^{0}(x)= \g(x), \quad x \in \O ,
\end{cases}
\end{equation*}
and the resolvent $J_{\tauhm\An} := \pa{\Id+\tauhm\An}^{-1}$ is a single-valued non-expansive operator on $L^p(\O)$ since $\An$ is $m$-accretive~\cite{Kato67}.
In addition, problem~\eqref{discretisationb} is well-posed as we state now.
\begin{lem}
Let $\g \in L^p(\O)$. Suppose that $\underline{\tau} := \inf\limits_{h} \tauh > 0$ or $\sum\limits_{h=1}^{+ \infty} \tauh^{\max(2,p)} = + \infty$, then the iterates of ~\eqref{discretisationb}, starting from $\g$, have a unique weak accumulation point $u^{\star} \in (\A)^{-1}(0)$. Moreover, if $\underline{\tau} > 0$, then for $h \geq 1$
\[
\norm{\An(\h)}_{L^p(\O)} \leq \frac{2\norm{\g-u^\star}_{L^p(\O)}}{(\underline{\tau}C_p)^{1/\max(p,2)} h^{1/\max(p,2)}} .
\]
\label{lem:existPdb}
\end{lem}
\bpf{}
$\An$ is accretive on $L^p(\O)$ (see the proof of~\cite[Theorem~6.7]{vaillo}). Moreover, it is well-known that for $p \in ]1,+\infty[$, $L^p(\O)$ is a uniformly convex and a uniformly smooth Banach space, whose convexity modulus verifies
\[
\delta_{L^p(\O)}(\varepsilon) \geq 
\begin{cases}
p^{-1}2^{-p}\varepsilon^{p}	& p \in [2,+\infty[ , \\
(p-1)\varepsilon^2/8		& p \in ]1,2] .
\end{cases}
\]
Thus, we are in position to apply~\cite[Theorem~3]{reich} to get uniqueness of the weak accumulation point. 

Let us turn to the rate. By $m$-accretiveness $\An$, $J_{\tauhm\An}$ is a single-valued operator on the entire $L^p(\O)$, and verifies for any $v,w \in L^p(\O)$ and $\lambda \in [0,1]$,
\begin{equation}
\label{eq:accretresolv}
\norm{J_{\tauhm\An}(v) - J_{\tauhm\An}(w)}_{L^p(\O)} \leq \norm{\lambda(v - w)+(1-\lambda)(J_{\tauhm\An}(v)-J_{\tauhm\An}(w))}_{L^p(\O)}.
\end{equation}
We now evaluate \eqref{eq:accretresolv} at $v=\H$, $w=u^\star$ and $\lambda=1/2$, and combine it with \cite[Corollary~2]{xu91}. This leads us to consider two possible cases.
\begin{enumerate}[label=(\alph*)]
\item $p \in ]2,+\infty[$: since $\h=J_{\tauhm\An}(\H)$ and $u^{\star}$ is a fixed point of $J_{\tauhm\An}$, and in view of \cite[Corollary~2, (3.4)]{xu91}, we have
\begin{align*}
\norm{\h - u^{\star}}_{L^p(\O)}^p &\leq \norm{\tfrac{1}{2}(\H-u^{\star}) + \tfrac{1}{2}(\h-u^{\star})}_{L^p(\O)}^p \\
				  &\leq \tfrac{1}{2}\norm{\H-u^{\star}}_{L^p(\O)}^p + \tfrac{1}{2}\norm{\h-u^{\star}}_{L^p(\O)}^p - 2^{-p}c_p\norm{\H-\h}_{L^p(\O)}^p \\
				  &\leq \norm{\H-u^{\star}}_{L^p(\O)}^p - 2^{-p}c_p\norm{\h-\H}_{L^p(\O)}^p,
\end{align*}
where we used non-expansiveness of $J_{\tauhm\An}$ to get the last inequality. $c_p = (1+\nu_p^{p-1})(1+\nu_p)^{1-p}$, where $\nu_p$ is the unique solution to $(p-2)\nu^{p-1}+(p-1)\nu^{p-2} = 1$, for $\nu \in ]0,1[$. Summing up these inequalities and using the fact that
\[
\norm{u_n^{h+1}-\h}_{L^p(\O)} \leq \norm{\h-\H}_{L^p(\O)}
\]
again by non-expansiveness of $J_{\tauhm\An}$, we arrive at
\[
\underline{\tau}h\norm{\An(\h)}_{L^p(\O)}^p \leq h\norm{\h-\H}_{L^p(\O)}^p \leq \sum_{i=1}^{h} \norm{u_n^{i}-u_n^{i-1}}_{L^p(\O)}^p \leq 2^p\norm{\g-u^{\star}}_{L^p(\O)}^p/c_p .
\]

\item $p \in ]1,2]$: using now \cite[Corollary~2, (3.7)]{xu91} and similar arguments to the first case, we get the inequality
\begin{align*}
\norm{\h - u^{\star}}_{L^p(\O)}^2 \leq \norm{\H-u^{\star}}_{L^p(\O)}^2 - 2^{-2}(p-1)\norm{\h-\H}_{L^p(\O)}^2 .
\end{align*}
Summing up again we end up with
\[
\underline{\tau}h\norm{\An(\h)}_{L^p(\O)}^2 \leq h\norm{\h-\H}_{L^p(\O)}^2 \leq \sum_{i=1}^{h} \norm{u_n^{i}-u_n^{i-1}}_{L^p(\O)}^2 \leq 4\norm{\g-u^{\star}}_{L^p(\O)}^2/(p-1) .
\]
\end{enumerate}
\epf
\begin{rem}
Observe that the assumption on the initial condition in Lemma~\ref{lem:existPdb} is weaker than that of Lemma~\ref{lem:existPdf} for $p \in ]1,2[$. As expected, the stability constraint needed on the time-step sequence is less restrictive than for the explicit/forward discretization.

\end{rem}
\begin{rem}
\begin{enumerate}[label=(\alph*)]
\item Observe that the assumption on the initial condition in Lemma~\ref{lem:existPdb} is weaker than that of Lemma~\ref{lem:existPdf}. 
\item As expected, the stability constraint needed on the time-step sequence is less restrictive than for the explicit/forward discretization.
\item Given that $\acc{\norm{u_n^{h+1}-\h}_{L^p(\O)}^p}_h$ is a decreasing and summable sequence, one can show that the rate $\norm{\An(\h)}_{L^p(\O)}=O(h^{-1/\max(p,2)})$ is in fact $\norm{\An(\h)}_{L^p(\O)} = o(h^{-1/\max(p,2)})$.
\end{enumerate}
\end{rem}

Equipped with this result, the proof of an analogue to Theorem~\ref{totally} in the implicit case is similar to that of the explicit case modulo the following change
\begin{equation*}
\M (x,t) = \sum_{h=1}^N \h(x) \chi_{]t_{h-1},t_h]} (t).
\end{equation*}

\subsection{Relation to Kobayashi type estimates} 
Consider the evolution problem
\begin{equation}\tag{CP}
\begin{cases}
u_t + A(t) u(t) \ni f(t), \\
u(0) = g.
\end{cases}
\label{cauchy}
\end{equation}
A problem of the form~\eqref{cauchy} is called an abstract Cauchy problem. The evolution problem \eqref{neumann} we deal with can be viewed as a particular case of \eqref{cauchy} in its autonomous-homogenous case, i.e. the operator $A(t) \equiv \A$ does not depend on time and $f \equiv 0$.

Problem~\eqref{cauchy} in the autonomous-homogenous \newcor{case} was studied by Kobayashi in~\cite{kobayashi}, where he constructed sequences of approximate solutions which converge in an appropriate sense to a solution to the differential inclusion. He provided an inequality that estimates the distance between arbitrary points of two independent sequences generated by the so called proximal iterations, from which, he derived quantitative estimates to compare the continuous and discrete trajectories using the backward Euler scheme. These estimates have similar flavour to ours when $K=K_n$. Later on, these results were generalized to the non-autonomous case as well as to the case where the trajectories are defined by two differential inclusions systems (i.e. different operators $A$); see \cite{alvarez} and references therein for a thorough review. The latter bounds, expressed in our notation, are provided only in terms of $\norm{\A(v)-\An(v)}_{L^p(\O)}$. We go further by exploiting the properties of our operators to get sharp estimates in terms of the data $\norm{\K-\Kn}_{L^p(\O^2)}$. This is more meaningful in our context where we recall that the goal is to study the fully discretized nonlocal $p$-Laplacian problem on graphs.

\section{Application to graph sequences}
\label{applications}
\subsection{Networks on simple graphs}
\label{simplegraphs}

A fully discrete counterpart of \eqref{neumann} on $\{G_n\}_n$ is then given by
\begin{equation}\tag{\textrm{$\mathcal{P}^{s,d}_n$}}
\begin{cases}
\displaystyle{\frac{\hi - \Hi}{\tauhm}} = \frac{1}{n}\sum\limits_{j:(i,j) \in E(G_n)} \abs{u_{j}^{h-1} - \Hi}^{p-2} (u_{j}^{h-1} - \Hi), & (i,h) \in [n] \times [N],\\
\\ 
u_i(0) = g_i^{0} , \quad i \in [n],
\end{cases}
\label{neumanndisc}
\end{equation} 
where 
\[
g_i^{0} = n \int_{\O_i^{(n)}} \g(x) dx 
\]
is the average value of $\g(x)$ on $\O_i^{(n)}$.

Let us recall that our main goal is to compare the solutions of the discrete and continuous models and establish some consistency results.  Since the two solutions do not live on the same spaces, it is practical to represent some intermediate model that is the continuous extension of the discrete problem, using the vector  $U^{h} = (u_1^{h}, u_2^{h} ,\cdots, u_n^{h})^T$  whose components uniquely solve the previous system~\eqref{neumanndisc} (see Lemma~\ref{lem:existPdf}) to obtain the following piecewise time linear interpolation on $\O \times [0,T]$ 
\begin{equation}
\U(x,t) = \frac{t_h - t }{\tauhm} \Hi + \frac{t-t_{h-1}}{\tauhm} \hi \quad  \text{if} \quad x \in \O_i^{(n)},\quad t\in ]t_{h-1},t_h],
\label{stepf}
\end{equation}
and the following piecewise constant approximation 
\begin{equation}
\M (x,t) = \sum_{i=1}^n\sum_{h=1}^N \Hi \chi_{]t_{h-1},t_h]} (t) \chi_{\O_i^{(n)}}(x).
\label{piecef}
\end{equation}
So that $\U(x,t)$ uniquely solves the following problem
\begin{equation}\tag{\textrm{$\mathcal{P}^s_n$}}
\begin{cases}
\frac{\partial}{\partial t } \U(x,t) = -\boldsymbol{\Delta}^{\C}_p(\M(x,t)) , & (x,t) \in \O \times ]0,T],\\
\U^{0}(x) = \g(x), \quad x \in \O ,
\end{cases}
\label{disc}
\end{equation}
where 
\[
\g(x) = g_i := n \int_{\O_i^{(n)}} g_n (x) dx \quad  \text{if} \quad x \in \O_i^{(n)} , i \in [n],
\]
$g_n(\cdot)$ being the initial condition taken in problem~\eqref{general} and $\C(x,y)$ is the piecewise constant function such that for $(x,y) \in  \O_{ij}^{(n)}$, $(i,j) \in [n]^2$
\begin{equation*}
\begin{cases}
n^2 \displaystyle{\int_{\O_{ij}^{(n)}}} K(x,y) dx dy  \quad \text{if} \quad \O_i^{(n)} \times  \O_j^{(n)} \cap \cl{\supp(K)} \ne \emptyset, \\
0 \qquad \text{otherwise}.
\end{cases}
\end{equation*}
As $G_n $ is a simple graph, $\C (\cdot, \cdot)$ is also a $\{0,1\}$-valued mapping.

\newcor{By analogy of what was done in \cite{medv}}, the rate of convergence of the solution of the discrete problem to the solution of the limiting problem depends on the regularity of the boundary $\bd{\cl{\supp(K)}}$ of the support closure. \newcor{Following ~\cite{medv}}, we recall the upper box-counting (or Minkowski-Bouligand) dimension of $\bd{\cl{\supp(K)}}$ as a subset of $\mathbb{R}^2$: 
\[
\rho := \dim_B(\bd{\cl{\supp(K)}}) = \limsup_{\delta \rightarrow 0} \frac{\log N_{\delta}(\bd{\cl{\supp(K)}})}{- \log \delta},
\]
where $ N_{\delta}(\bd{\cl{\supp(K)}}) $ is the number of cells of a $( \delta \times \delta )$-mesh that intersect $ \bd{\cl{\supp(K)}} $ 
(see \cite{falconer}).
\begin{cor}
Suppose that $p \in ]1,+\infty[$, $g\in L^{\infty} (\O)$, and 
\[
\rho \in [0,2[.
\]
Let $u$ and $\U$ denote the  functions corresponding to the solutions of  \eqref{neumann} and \eqref{disc}, respectively.

Then for any $\epsilon >0$ there exists $N(\epsilon) \in \mathbb{N}$ such that for any $n \ge N(\epsilon)$
\begin{equation}
\norm{u- \U}_{C(0,T;L^p(\O))} \le C \left (  \norm{g - g_n}_{L^p(\O)} + n^{-((2-\rho)/p - \epsilon)} \right ) + O(\tau),
\end{equation}
where the positive constant $C$ is independent of $n$.
\label{theo1}
\end{cor}
\bpf{}
By Theorem~\ref{totally}, we have
\begin{equation}
\norm{u- \U}_{C(0,T;L^p(\O))} \le C \left (  \norm{g - g_n}_{L^p(\O)} +\norm{ g_n-\g}_{L^p(\O)} + \norm{K - \C}_{L^p(\O)}  \right )+ O(\tau).
\label{simple}
\end{equation}
Since both~\eqref{disc} and~\eqref{neumanndisc} problems share the same initial data, we have that $ \norm{g_n-\g}_{L^p(\O)} = 0$. It remains to estimate $ \norm{ K - \C  }_{L^{p}(\O)}$. \newcor{To do this, we follow the same proof strategy as in \cite[Theorem 4.1]{medv} }. For that, consider the set of discrete cells $\O_{ij}^{(n)}$ overlying the boundary of the support of $K$ 
\[
S(n) = \ens{(i,j) \in [n]^2}{\O_{ij}^{(n)} \cap \bd{\cl{\supp(K)}} \ne \emptyset} \qandq C(n) = \abs{S(n)}.
\]
For any $\epsilon > 0$ and sufficiently large $n$, we have 
\[
C(n) \le n^{\rho + \epsilon}.
\]
It is easy to see that $K$ and $\C$ coincide almost everywhere on cells $\O_{ij}^{(n)}$ for which $(i,j)  \notin S(n)$. Thus for any $\epsilon > 0$ and all sufficiently large $n$, we have 
\begin{equation}
 \norm{ K- \C }_{L^{p}(\O^2)}^p = \int_{\O^2} |K(x,y)- \C(x,y)|^p dx dy \le C(n) n^{-2} \le n^{-(2-\rho-\epsilon)}.
 \label{sup}
\end{equation}
Assembling \eqref{simple} and \eqref{sup}, the desired result holds.
\epf

\subsection{Networks on weighted graphs}
\label{weightedgraphs}
\subsubsection{Networks on  $K / \mathcal{Q}_n$}
We consider the totally discrete counterpart of \eqref{neumann} on $K / \mathcal{Q}_n$
\begin{equation}\tag{\textrm{$\hat{\mathcal{P}}^{w,d}_{n}$}}
\begin{cases}
\displaystyle{\frac{\hi - \Hi}{\tauhm}}= \frac{1}{n} \sum\limits_{j=1}^{n} (\k)_{ij}  \abs{u_j^{h-1} - \Hi}^{p-2} (u_j^{h-1} - \Hi), & (i,h) \in [n] \times [N],
\\ 
u_i (0) = g_i^{0} , \quad i \in [n],
\end{cases}
\label{pond}
\end{equation}
where $\k$ is defined in \eqref{eq:kweight} and $g_i^{0}$ is the average value of $\g(x)$ on ${\O^{(n)}_{i}} $.

Combining the piecewise constant function $\U$  in~\eqref{stepf} with $\M$ in~\eqref{piecef}, we rewrite \eqref{pond} as 
\begin{equation}\tag{\textrm{$\hat{\mathcal{P}}^{w}_n$}}
\begin{cases}
\frac{\partial}{\partial t } \U(x,t) = -\boldsymbol{\Delta}^{\D}_p(\M(x,t)), & (x,t) \in \O \times ]0,T],\\
\U^{0}(x) = \g(x), \quad x \in \O ,
\label{pondd}
 \end{cases}
 \end{equation}
 
 where $\D$ and $\g$ are the piecewise constant functions such that 
 \[
 \D(x,y) = (\k)_{ij} \quad \text{for} \quad (x,y) \in \O^{(n)}_{i} \times \O^{(n)}_{j},
 \]
 \[
 \g(x) = g_{i} \quad \text{for} \quad x \in \O^{(n)}_{i} , i \in [n].
 \]

\newcor{As already emphasized in \cite[Remark~5.1]{medv}, it is instructive to note that~\eqref{pondd} can be viewed as the time discretized Galerkin approximation of problem~\eqref{neumann}.}


\begin{cor}
Suppose that $p \in ]1,+\infty[$, $K: \O^2 \to [0,1]$ is a symmetric measurable function, and $g \in L^{\infty}(\O)$. Let $u$ and $\U$ be the solutions of~\eqref{neumann} and~\eqref{pondd}, respectively.
Then 
\begin{equation}
\norm{u-\U}_{C(0,T;L^p(\O))} \underset{n \to \infty, \tau \to 0}{\longrightarrow} 0 .
\label{asympto}
\end{equation}
\end{cor}
\bpf{} \newcor{This proof strategy was used in \cite[Theorem 5.2]{medv}}. For fixed $(i,j) \in [n]^2$, it is easy to see that $\{\O_{ij}^{(n)}\}_{n}$ is a decreasing sequence, $\bigcap\limits_{n=1}^{\infty} \O_{ij}^{(n)} = \{(x,y)\}$, and 
\[
 (\k)_{ij} = \frac{1}{\abs{\O_{ij}^{(n)}}} \int_{\O_{ij}^{(n)}}  K_{n} (x,y) dx dy .
\]
Then, by the Lebesgue differentiation theorem (see e.g. \cite[Theorem~3.4.4]{pardoux}), we have 
\[
\D \underset{n \to \infty}{\longrightarrow} K, 
\]
almost everywhere on $\O^2$, whence, using the same arguments on $\R$, we have also that $g_n  \underset{n \to \infty}{\longrightarrow} g$ almost everywhere on $\O$. Thus, combining Theorem~\ref{totally} and statement~\ref{scheffe} in Theorem~\ref{main}, the desired result follows.
\epf

To quantify the rate of convergence in  \eqref{asympto}, we need to add some supplementary assumptions on the kernel $K$ and the initial data $g$. To do this, we introduce the Lipschitz spaces $\Lip(s,L^p(\O^d))$, for $d \in \{1,2\}$, which contain functions with, roughly speaking, $s$ "derivatives" in $L^p(\O^d)$~\cite[Ch.~2, Section~9]{devorelorentz93}. 
\begin{defi}\label{def:lipspaces}
For $F \in L^p(\O^d)$, $p \in [1,+\infty]$, we  define the (first-order) $L^p(\O^d)$ modulus of smoothness by
\begin{equation}
\omega(F,h)_p := \sup_{\bs z \in \R^d, |\bs z| < h} \pa{\int_{\bs x,\bs x +\bs z \in \O^d}\abs{F(\bs x + \bs z)-F(\bs x)}^pd\bs x}^{1/p} .
\label{modsmooth}
\end{equation}
The Lipschitz spaces $\Lip(s,L^p(\O^d))$ consist of all functions $F$ for which
\[
\abs{F}_{\Lip(s,L^p(\O^d))} := \sup_{h > 0} h^{-s} \omega(F,h)_p < +\infty .
\]
\end{defi}
We restrict ourselves to values $s \in ]0,1]$ as for $s > 1$, only constant functions are in $\Lip(s,L^p(\O^d))$. It is easy to see that $\abs{F}_{\Lip(s,L^p(\O^d))}$ is a semi-norm. $\Lip(s,L^p(\O^d))$ is endowed with the norm
\[
\norm{F}_{\Lip(s,L^p(\O^2))} :=  \norm{F}_{L^p(\O^2)} +  \abs{F}_{\Lip(s,L^p(\O^d))} .
\]
The space $\Lip(s,L^p(\O^2))$ is the Besov space $\mathbf{B}^s_{p,\infty}$~\cite[Ch.~2, Section~10]{devorelorentz93} which are very popular in approximation theory. In particular, $\Lip(1,L^1(\O^d))$ contains the space $\BV(\O^d)$ of functions of bounded variation on $\O^d$, i.e. the set of functions $F \in L^1(\O^d)$ such that their variation is finite:
\begin{equation*}
V_{\O^2}(F) := \sup_{h > 0}h^{-1}\sum_{i=1}^d\int_{\O^d}\abs{F(\bs x + he_i)-F(\bs x)}d\bs x < + \infty 
\end{equation*}
where $e_i, i \in \{1,d\}$ are the coordinate vectors in $\R^d$; see~\cite[Ch.~2, Lemma~9.2]{devorelorentz93}. Thus Lipschitz spaces are rich enough to contain functions with both discontinuities and fractal structure.

Let us define the piecewise constant approximation of a function $F \in L^p(\O^2)$ (a similar reasoning holds on $\O$),
\[
\hat{F}_n(x,y) := \frac{1}{\abs{\O_{ij}^{(n)}}}\sum_{ij} \pa{\int_{\O^2}  F(x',y')\chi_{\O_{ij}^{(n)}} (x',y') dx' dy'} \chi_{\O_{ij}^{(n)}}(x,y),
\]
where $\chi_{\O_{ij}^{(n)}}$ is the characteristic function of $\O_{ij}^{(n)}$. Clearly, $\hat{F}_n$ is nothing but the projection $\mathbf{P}_{V_{n^2}}(F)$ of $F$ on the $n^2$-dimensional subspace $V_{n^2}$ of $L^{p}(\O^2)$ defined as $V_{n^2} = \Span \ens{\chi_{\O_{ij}^{(n)}}}{(i,j)\in [n]^2}$.


\begin{lem} 
There exists a constant $C$ such that for all $F \in \Lip(s,L^p(\O^2))$, $s \in ]0,1]$, $p \in [1,+\infty]$,\begin{equation}
\norm{F - \hat{F}_n}_{L^p(\O^2)} \le C \frac{\abs{F}_{\Lip(s,L^p(\O^2))}}{n^s}.
\label{eq:lipspaceapprox}
\end{equation}
In particular, if $F \in \BV(\O^2) \cap L^\infty(\O^2)$, then
\begin{equation}
\norm{F - \hat{F}_n}_{L^p(\O^2)} \le \frac{\Bpa{C\bpa{2\norm{F}_{L^{\infty}(\O^2)}}^{p-1} V_{\O^2}(F)}^{1/p}}{n^{1/p}}.
\label{eq:bvapprox}
\end{equation}
\label{lem:spaceapprox}
\end{lem}
Similar bounds hold for $g$.

\bpf{} 
Using the general bound \cite[Ch.~7, Theorem~7.3]{devorelorentz93} for the error in spline approximation, and in view of Definition~\ref{def:lipspaces}, we have
\[
\norm{F - \hat{F}_n}_{L^p(\O^2)} \leq C \omega(F,1/n)_p = C n^{-s} (n^s \omega(F,1/n)_p) \leq C n^{-s} \abs{F}_{\Lip(s,L^p(\O^2))} .
\]
As for \eqref{eq:bvapprox}, we know that $\BV(\O^2) \subset \Lip(1,L^1(\O^2))$. Thus, combining Lemma~\ref{inclusion}, \eqref{eq:lipspaceapprox} and \cite[Ch.~2, Lemma~9.2]{devorelorentz93}, we get 
\begin{align*}
\norm{F - \hat{F}_n}_{L^p(\O^2)} 
&\leq \norm{F - \hat{F}_n}_{L^{\infty}(\O^2)}^{1-\frac{1}{p}} \norm{F - \hat{F}_n}_{L^1(\O^2)}^{\frac{1}{p}} \\
&\leq \bpa{2\norm{F}_{L^{\infty}(\O^2)}}^{1-\frac{1}{p}} (C V_{\O^2}(F))^{1/p} n^{-1/p} .
\end{align*}
\epf
The second claim \eqref{eq:bvapprox} can also be proved using \cite[Lemma~3.2(3)]{Lucierlectures96}.\\


We are now in position to state the following error bound.
\begin{cor}
Suppose that $p \in ]1,+\infty[$, $K: \O^2 \to [0,1]$ is a symmetric and measurable function in $\Lip(s,L^p(\O^2))$, and $g \in \Lip(s,L^p(\O)) \cap L^\infty(\O)$, $s \in ]0,1]$. Let $u$ and $\U$ be the solutions of \eqref{neumann} and \eqref{pondd} respectively. Then 
\begin{equation}
\norm{u - \U}_{C(0,T;L^p(\O))} \le  O(n^{-s})+ O(\tau).
\end{equation}
If $\Lip(s,L^p(\O^2))$ is replaced with $\BV(\O^2)$, then the rate becomes 
\begin{equation}
\norm{u - \U}_{C(0,T;L^p(\O))} \le  O(n^{-1/p})+ O(\tau).
\end{equation}
\label{weight}
\end{cor} 
\bpf{}By Theorem~\ref{totally}, we have
\begin{equation*}
\norm{u- \U}_{C(0,T;L^p(\O))} \le C \left (  \norm{g - g_n}_{L^p(\O)} +  \norm{g_n-\g}_{L^p(\O)} + \norm{K - \D}_{L^p(\O)}  \right ) + O(\tau).
\end{equation*}
Since the initial conditions for both~\eqref{pond} and~\eqref{pondd} stem from the same initial data, we have that $ \norm{g_n-\g}_{L^p(\O)} = 0$. The claimed rates then follow from Lemma~\ref{lem:spaceapprox} since $\D=\mathbf{P}_{V_{n^2}}(K)$ and $g_n = \mathbf{P}_{V_{n}}(g)$.
\epf

\subsubsection{The limit as $p \to \infty$}
Let us consider the numerical fully discrete approximation of the problem~\eqref{neumann} using the function $\k$ defined in~\eqref{eq:kweight}
\begin{equation}
\begin{cases}
\displaystyle{\frac{ U_{i,h}^p - U_{i,h-1}^p }{\tauhm}} = \frac{1}{n} \sum\limits_{j=1}^{n} (\k)_{ij}  \abs{U_{j,h-1}^p  - U_{i,h-1}^p }^{p-2} (U_{j,h-1}^p  - U_{i,h-1}^p ), & (i,h) \in [n] \times [N],\\
U^p_{i,0} = g_i^{0}, \quad i \in [n],
\end{cases}
\label{pproblem}
\end{equation}
where the vector $U^p \in \R^{nN}$. This problem is associated to the energy functional 
\[
F_{p}(V) = \frac{1}{2pn^2} \sum\limits_{i= 1}^{n}\sum\limits_{j=1}^{n} (\k)_{ij} \abs{V_j- V_i}^{p},
\]
in the Euclidean space $H := \R^{n}$.

As before, we consider the linear interpolation of $U^p$ as follows 
\begin{equation}
\R^n \ni \check{U}^p(t)= \frac{t_{h} - t}{\tauhm} U_{h-1}^p + \frac{t-t_{h-1}}{\tauhm}U_{h}^p ,\quad t \in ]t_{h-1},t_h] ,
\label{convexx}
\end{equation}
and a piecewise constant approximation
\begin{equation}
\R^n \ni \bar{U}^p(t) =  U_h^p,\quad t \in ]t_{h-1},t_h].
\label{piecewisee}
\end{equation}
Consequently, $\check{U}^p(.)$ obeys the following evolution equation
\begin{equation}
\begin{cases}
\displaystyle{\frac{d \check{U}^p(t) }{dt}} = \frac{1}{n} \sum\limits_{j=1}^{n} (\k)_{ij}   \abs{\bar{U}_j^p(t)  - \bar{U}_i^p(t) }^{p-2} (\bar{U}_j^p(t)  - \bar{U}_i^p(t) ), & (i,t) \in [n] \times ]0,T] ,\\
U_i^p(0) = g_i^{0}, \quad i \in [n].
\end{cases}
\end{equation}
Now we define 
\begin{equation}
\begin{cases}
\displaystyle{\frac{d U^p(t) }{dt}} = \frac{1}{n} \sum\limits_{j=1}^{n} (K_n)_{ij}  \abs{U_j^p(t)  - U_i^p(t) }^{p-2} (U_j^p(t)  - U_{i}^p(t)), & (i,t) \in [n] \times ]0,T] ,\\
U_i^p(0) = g_i^{0}, \quad i \in [n].
\end{cases}
\end{equation}

To avoid triviality, we suppose that $\supp(\k)  \ne \emptyset$, and define the non-empty compact convex set
\[
S_{\infty} = \ens{v \in \R^{n N}}{\abs{v_j - v_i} \le 1, \quad \text{for} \quad (i,j) \in \supp(\k)},
\]
where the subscript $\infty$ will be made clear shortly. Indeed, taking the limit as $p\to \infty$ of $F_p$, one clearly sees that this limit is $\imath_{S_{\infty}}$, where the latter is the indicator function of $S_{\infty}$, i.e.
\begin{equation*}
\imath_{S_{\infty}}(v) = 
\begin{cases}
0\quad \text{if} \quad v \in S_{\infty},\\
+ \infty \quad \text{otherwise}.
\end{cases}
\end{equation*}
Then, the nonlocal time continuous limit problem can be written as 
\begin{equation}\tag{\textrm{$\mathcal{P}^{\infty}$}}
\begin{cases}
\displaystyle{\frac{d U^{\infty} }{dt}}+ \text{N}_{S_{\infty}}(U^{\infty}(t)) \ni 0, & t \in ]0,T] ,\\
U_i^{\infty} (0) = g_i^{0}, \quad i \in [n],
\end{cases} 
\label{infinity}
\end{equation}
where  $\text{N}_{S_{\infty}}$ denotes the normal cone of $S_{\infty}$, defined by 
\[
\text{N}_{S_{\infty}}(v) = 
\begin{cases}
\ens{\eta \in H}{\pds{\eta}{w - v} \le 0 , \forall w \in H} & \text{if } v \in S_{\infty}, \\
\emptyset & \text{otherwise},
\end{cases}
\]
where $\pds{\cdot}{\cdot}$ denotes the inner product on the Hilbert space $H$.
\begin{theo} Suppose that $\supp(\k)  \ne \emptyset$ and $g^0 \in S_{\infty}$. Let $\check{U}^p$ be the solution of~\eqref{pproblem}. If $U^{\infty}$ is the unique solution to ~\eqref{infinity}, then
\begin{equation}
\lim\limits_{p \to \infty}\lim\limits_{\tau \to 0}\sup_{t \in [0, T]} \abs{\check{U}^p(t)- U^{\infty}(t)} = 0,
\label{limitt}
\end{equation}
\newcor{where $\tau = \max\limits_{h \in [N]}\tauh$ is is the maximal size of intervals in the partition of $[0,T]$.}
\label{theolimit}
\end{theo}
\begin{rem} Before carrying out the proof of Theorem~\ref{theolimit}, note that one cannot interchange the order of limits; the limit as $\tau \to 0$ must be taken before the limit as $p \to \infty$. The reason will be clarified in the proof.
\end{rem}
\bpf{}Using the triangle inequality, we have 
\[
\abs{\check{U}^p(t)- U^{\infty}(t)} \le \abs{\check{U}^p(t) - U^p(t)} + \abs{U^p (t)- U^{\infty}(t)}.
\]
First, proceeding exactly as in the proof of Theorem~\ref{totally}, and more precisely inequality~\eqref{approx}, we get 
\begin{equation}
\abs{\check{U}^p(t)- U^p(t)} \leq C' \tau
\label{approxi}
\end{equation}
for $C' \geq 0$. Since the constant $C'$ in~\eqref{approxi} depends on $p$, we first take the limit as $\tau \to 0$, to get 
\begin{equation}
\lim\limits_{\tau \to 0}\sup_{t \in [0, T]} \abs{\check{U}^p(t)- U^p(t)} = 0
\label{limit1}
\end{equation}
Now, arguing as in \cite[Theorem~3.2]{rossi} (which in turn relies on \cite[Theorem~3.1]{brezis}), we have additionally that 
\begin{equation}
\lim\limits_{p \to \infty} \sup_{t \in [0, T]} \abs{U^p(t)- U^{\infty}(t)} = 0.
\label{limit2}
\end{equation}
Hence, the combination of~\eqref{limit1} and~\eqref{limit2} yields~\eqref{limitt}.
\epf

\begin{rem}
Note that we get the same result when dealing with the implicit Euler scheme, following the changes mentioned in Section~\ref{subsec:Pdb}.
\end{rem}
\subsubsection{Networks on $\mathbb{G} (X_{n},K)$}
The analysis of the problem \eqref{neumann} on  $\mathbb{G} (X_{n},K)$ remains the same modulo the definition of the piecewise constant approximation 
\[
\breve{K}_n^w(x,y) = (\breve{K_{n}})_{ij} \quad \text{for} \quad (x,y) \in \O_{ij}^{(n)},
\]
where we recall $\breve{K_{n}}$ from~\eqref{collocation}.
The fully discrete counterpart of~\eqref{neumann} on  $\mathbb{G} (X_{n},K)$  is given by 
\begin{equation}\tag{\textrm{$\breve{\mathcal{P}}^{w,d}_{n}$}}
\begin{cases}
\displaystyle{\frac{\hi - \Hi}{\tau}}= \frac{1}{n} \sum\limits_{j=1}^{n}   (\breve{K_{n}})_{ij}  \abs{\hi - \Hi}^{p-2} (u_j^{h-1} - \Hi), & (i,h) \in [n] \times [N] ,
\\ 
u_i (0) = g_i^{0}, \quad i \in [n].
\end{cases}
\label{pondGX}
\end{equation}
It is worth mentioning that~\eqref{pondGX} is the time discretized approximation of the problem~\eqref{neumann} using the collocation method. Roughly speaking, it is about the projection of~\eqref{neumann} on $X_n$ (see~\eqref{collocation}) via the interpolation operator  $P_n: L^{\infty}(\O) \to X_n$ which to each $u(t_h, .) \in L^{\infty}(\O)$ associates the unique function $f(t_h, .)$ such that for all $i \in [n]$, $u(t_h ,\frac{i}{n})  = f(t_h ,\frac{i}{n})$. See \cite{gabriel} for more details.

We assume further that the kernel $K$ is almost everywhere continuous on $\O^2$. By construction of $\breve{K}_n^w$ (see~\eqref{collocation}),
\[
\breve{K}_n^w(x,y) \to K(x,y), \quad \text{as} \quad n \to \infty,
\]
at every point of continuity of $K$, i.e., almost everywhere. Thus, using the Sheffe-Riesz theorem, we have 
\[
\norm{K-\breve{K}_n^w}_{L^p(\O^2)} \to 0 \quad \text{as} \quad n \to \infty.
\]
Thereby, the proof of Corollary~\ref{weight} applies to the situation at hand. Hence, we have the following result.
\begin{cor}
Suppose that $p \in ]1,+\infty[$, $K: \O^2 \to [0,1]$ is a symmetric measurable function, which is continuous almost everywhere on $\O^2$, and $g \in L^{\infty}(\O)$. Let $u$ be the solution of~\eqref{neumann}, and $\U$ be the piecewise constant extension as in \eqref{stepf} using the sequence in \eqref{pondGX}. Then
\[
\norm{u - \U}_{C(0,T;L^p(\O))} \to  0 \quad \text{as} \quad n \to \infty.
\]
\end{cor}
\begin{rem}
The result of Theorem~\ref{theolimit} remains the same for this graph model taking the kernel $(\breve{K_{n}})_{ij}$ instead of $(\hat{K_{n}})_{ij}$.
\end{rem}

\appendix
\section{Proof of Theorem~\ref{theo:exist}}
\label{subsec:prooftheoexist}
 
As stated above, the proof follows the lines of that of \cite[Theorem~6.7]{vaillo}. It relies on arguments from nonlinear semigroup theory (and in particular resolvents of accretive operators in Banach spaces). To apply the same arguments as for \cite[Theorem~6.7]{vaillo}, we need the following two lemmas that extend \cite[Lemmas~6.5 and 6.6]{vaillo} to the case of a symmetric, nonnegative and bounded kernel $K$.

\begin{lem} For every $u,v \in L^p(\O)$,
\begin{equation*}
\begin{split}
&- \int_{\O}\int_{\O} K(x,y) \abs{u(y) - u(x) }^{p-2} (u(y) - u(x) ) dy v(x) dx \\
&= \frac{1}{2} \int_{\O}\int_{\O} K(x,y) \abs{u(y) - u(x) }^{p-2} (u(y) - u(x) )(v(y) - v(x)) dy dx .
\end{split}
\end{equation*}
\label{monotone}
\end{lem}
From this lemma the following monotonicity result can be deduced.
\begin{lem} 
Let $T : \mathbb{R} \rightarrow \mathbb{R}$ be a nondecreasing function. Then 
\begin{enumerate}[label=(\roman*)]
\item  For every $u,v \in L^p(\O) $ such that $T(u-v) \in L^p(\O)$, we have 
\begin{equation}
\begin{split}
&\int_{\O} (\A u(x) - \A v(x)) T(u(x) - v(x) ) dx \\
&= \frac{1}{2} \int_{\O}\int_{\O} K(x,y) (T(u(y) - v(y)) - T(u(x) - v(x)))  \\
&\times \left (\abs{u(y) - u(x) }^{p-2} (u(y) - u(x) ) - \abs{v(y) - v(x) }^{p-2} (v(y) - v(x) )\right ) dy dx.
\end{split}
\label{lem1}
\end{equation}
\item Moreover, if $T$ is bounded  \eqref{lem1} holds for every $u,v \in \dom (\A)$.
\end{enumerate}
\label{monotone1}
\end{lem}
\subsection{Proof of Lemma~\ref{monotone}}

\bpf{}
Let $\O'$ be a bounded subset of $\mathbb{R}$ and let $\Gamma \subset \mathbb{R} \setminus \interop{ \O'}$.

For $\a : (\O' \cup \Gamma ) \times (\O' \cup \Gamma ) \rightarrow \R $ , $u : \O' \cup \Gamma   \rightarrow \R $, and $f: (\O' \cup \Gamma ) \times (\O' \cup \Gamma ) \rightarrow \R$. We define as in \cite{bord} the following generalized nonlocal operators 
\begin{enumerate}[label=(\alph*)]
\item \textbf{Generalized gradient}
\[
\mathcal{G}(u) (x,y) := (u(y) - u(x) ) \a(x,y), \quad x,y \in \O' \cup \Gamma,
\]
\item \textbf{Generalized nonlocal divergence}
\[
\mathcal{D}(f) (x,y) := \int_{\O' \cup \Gamma} (f(x,y) \a (x,y) - f(y,x)  \a (y,x)) dy, \quad x \in \O',
\]
\item \textbf{Generalized normal component}
\[
\mathcal{N}(f) (x,y) := - \int_{\O' \cup \Gamma} (f(x,y) \a (x,y) - f(y,x)  \a (y,x)) dy, \quad x \in \Gamma.
\]
\end{enumerate}
With the above notation in place, the authors in \cite{bord} prove that for $v : \O' \cup \Gamma \rightarrow \R $ and $s : \O' \cup \Gamma  \times \O' \cup \Gamma  \rightarrow \R$\newcor{,} the following identity holds 
\begin{equation}
\int_{\O'} v \mathcal{D}(s) dx + \int_{\O' \cup \Gamma } \int_{\O' \cup \Gamma } s \mathcal{G}(v) dy dx = \int_{\Gamma} v \mathcal{N} (s) dx.
\label{gen}
\end{equation}
Let $\mu :( \O' \cup \Gamma ) \times (\O' \cup \Gamma) \rightarrow \R $ be given by 
\[
\mu (x, y) := |\a(x,y)|^p.
\]
 In our particular case $\mu$ is the kernel $K(\cdot,\cdot)$, so that we suppose that $\a$ is symmetric.
 Hence, the following identity 
\[
\mathcal{D}(|\mathcal{G}(u)|^{p-2} \mathcal{G}(u)) = \mathcal{L}_p u := 2 \int_{\O' \cup \Gamma} |u(y) - u(x) |^{p-2} (u(y) - u(x)) \mu (x,y) dy
\]
 was also shown in \cite[(5.3)]{bord} for $p=2$. The general case was proved in \cite{general}, that is 
 \begin{equation}
 \mathcal{L}_p u = \mathcal{D}(\abs{\mathcal{G}(u)}^{p-2} \mathcal{G}(u)). 
 \label{gen1}
 \end{equation}
 The equality holds whenever both sides are finite.
 
Applying \eqref{gen} with $s(x,y) = \abs{\mathcal{G}(u)}^{p-2} \mathcal{G}(u)(x,y)$ and using the identity \eqref{gen1}, we obtain 

\begin{equation*}
\int_{\O'}  \mathcal{L}_p (u)  v dx + \int_{ \O' \cup \Gamma} \int_{ \O' \cup \Gamma} (\abs{\mathcal{G}(u)}^{p-2} \mathcal{G}(u)) . \mathcal{G}(v)dx dy = \int_{\Gamma}  \mathcal{N} ( \abs{\mathcal{G}(u)}^{p-2} \mathcal{G}(u))  v  dx.
\end{equation*}
Hence
\begin{equation*}
\begin{split}
  \int_{\O'} \mathcal{L}_p v dx &= - \int_{ \O' \cup \Gamma} \int_{ \O' \cup \Gamma} (\abs{\mathcal{G}(u)}^{p-2} \mathcal{G}(u))  \mathcal{G}(v)dx dy  +  \int_{\Gamma} v \mathcal{N} ( \abs{\mathcal{G}(u)}^{p-2} \\
 &= - \int_{ \O' \cup \Gamma} \int_{ \O' \cup \Gamma} (\abs{\mathcal{G}(u)}^{p-2} \mathcal{G}(u)) \mathcal{G}(v)dx dy \\
 & + \int_{\Gamma} \left( -  \int_{ \O' \cup \Gamma}   \abs{\mathcal{G}(u)}^{p-2} \mathcal{G}(u) (x,y)\a (x,y) -  \abs{\mathcal{G}(u)}^{p-2} \mathcal{G}(u)(y,x) \a(y,x) dy\right ) v dx\\
&= - \int_{ \O' \cup \Gamma} \int_{ \O' \cup \Gamma} \abs{\mathcal{G}(u)}^{p-2} \mathcal{G}(u)) \mathcal{G}(v)dx dy \\
&-  \int_{\Gamma} \int_{ \O' \cup \Gamma}  \a(x,y) \left (    \abs{\mathcal{G}(u)}^{p-2} \mathcal{G}(u) (x,y) -   \abs{\mathcal{G}(u)}^{p-2} \mathcal{G}(u) (y,x)\right ) dy  v dx \\
&= - \int_{ \O' \cup \Gamma} \int_{ \O' \cup \Gamma} \abs{\mathcal{G}(u)}^{p-2} \mathcal{G}(u)) \mathcal{G}(v)dx dy -  \int_{\Gamma} \mathcal{L}_p(u)  v dx.
 \end{split}
 \end{equation*}
 Thus 
 \begin{equation}
 \int_{\O' \cup \Gamma}  \mathcal{L}_p(u)  v dx  = -  \int_{ \O' \cup \Gamma} \int_{ \O' \cup \Gamma}|\mathcal{G}(u)|^{p-2} \mathcal{G}(u)) \mathcal{G}(v)dx dy.
 \label {gen2}
 \end{equation}
Replacing $\mathcal{G}$ with its form in \eqref{gen2} and taking $\O = \O' \cup \Gamma $ as this nonlocal integration formula does not contain any boundary terms, so that, the values of $u$ could be nonzero on the domain $\Gamma$ without affecting the formula, we get the desired result.
\epf
\subsection{Proof of Lemma~\ref{monotone1}}

\bpf{}
\begin{enumerate}[label=(\roman*)]
\item We have
\begin{align*}
& \int_{\O} (\A u(x) - \A v(x)) T(u(x) - v(x)) dx \\
&= \int_{\O} \left( - \int_{\O} K(x,y) \abs{u(y) - u(x)}^{p-2} (u(y) - u(x)) dy \right) T(u(x) - v(x) ) dx\\
&+ \int_{\O} \left (\int_{\O} K(x,y) \abs{v(y) - v(x)}^{p-2} (v(y) - v(x)) dy \right ) T(u(x) - v(x) ) dx \\
&= -  \int_{\O}  \int_{\O} K(x,y)   (\abs{u(y) - u(x)}^{p-2} (u(y) - u(x)) -\\
&  \abs{v(y) - v(x)}^{p-2} (v(y) - v(x))  ) dy T( u(x) - v(x) ) dx \\
&= -  \int_{\O}  \int_{\O} K(x,y) \abs{u(y) - u(x)}^{p-2} (u(y) - u(x)) dy T(u(x) - v(x) ) dx - \\
&-   \int_{\O}  \int_{\O} K(x,y) \abs{v(y) - v(x)}^{p-2} (v(y) - v(x)) dy  T(u(x) - v(x) )dx \\
&=  \frac{1}{2}  \int_{\O}  \int_{\O} K(x,y) \abs{u(y) - u(x)}^{p-2} (u(y) - u(x)) (T(u(y) - v(y) )-  T(u(x) - v(x) ) dx dy \\
&-  \frac{1}{2}  \int_{\O}  \int_{\O} K(x,y)  \abs{v(y) - v(x)}^{p-2} (v(y) - v(x)) (T(u(y) - v(y) )-  T(u(x) - v(x) )dx dy\\
&= -  \frac{1}{2}  \int_{\O}  \int_{\O}  K(x,y)( \abs{u(y) - u(x)}^{p-2} - |v(y) - v(x)|^{p-2} (v(y) - v(x))) \\
& \times  (T(u(y) - v(y) )-  T(u(x) - v(x) )dx dy.
\end{align*}
\item If $T$ is bounded , we have 
\[ 
\forall u,v \in \dom (\A) ,  \hspace{0.5cm}T(u-v) \in L^p(\O).
\]
\end{enumerate}
\epf

\section{Mean value theorem for continuous functions}
The following lemma states a generalization of the Lagrange mean value theorem retaining only the continuity assumption, but weakening the differentiability hypothesis.

\begin{lem}\label{lem:genmvt}
Suppose that the real-valued function $f$ is continuous on $[a,b]$, where $a < b$, both $a$ and $b$ being finite. If the right and left-derivatives $f'_{+}$ and $f'_{-}$ exist as extended-valued functions on $]a,b[$, then there exists $c \in ]a,b[$ such that either
\[
f'_{+}(c) \leq \frac{f(b)-f(a)}{b-a} \leq f'_{-}(c)
\]
or
\[
f'_{-}(c) \leq \frac{f(b)-f(a)}{b-a} \leq f'_{+}(c) .
\]
If moreover $f'_{+}$ and $f'_{-}$ coincide on $]a,b[$, then $f$ is differentiable at $c$ and
\[
f(b)-f(a)=f'(c)(b-a) .
\]
\end{lem}
\bpf{}
From \cite[p.~115]{DiazVyborny64} (see also \cite{Young1909}), we have under the sole continuity assumption of $f$ on $[a,b]$ that either
\[
\frac{f(c+h)-f(c)}{h} \leq \frac{f(b)-f(a)}{b-a} \leq \frac{f(c)-f(c-d)}{d} 
\]
or
\[
\frac{f(c)-f(c-d)}{d} \leq \frac{f(b)-f(a)}{b-a} \leq \frac{f(c+h)-f(c)}{h},
\]
for all $h > 0$ and $d > 0$ such that $(c+h,c-d) \in ]a, b[^2$. Passing to the limit as $h \to 0^+$ and $d \to 0^+$ (the limits exist in $[-\infty,+\infty]$ by assumption), we get our inequalities. When $f'_{+}$ and $f'_{-}$ coincide on $]a,b[$, and in particular at $c$, the inequalities become an equality $f'_{+}(c)=f'_{-}(c)=\frac{f(b)-f(a)}{b-a}$, and the derivative at $c$ is finite, whence differentiability follows.
\epf

Let us apply this result to $f: t \in \R \mapsto \abs{t}^{p-2}t$, $p > 1$. $f$ is a continuous\footnote{Observe that $f$ is not even continuous at $0$ when $p=1$, and thus Lemma~\ref{lem:genmvt} cannot be applied when $0 \in [a,b]$.} monotonically increasing and odd function on $\R$ . It is moreover everywhere differentiable for $p \geq 2$, and for $p \in ]1,2[$ it is differentiable except at $0$, where $f'_{+}(0)=f'_{-}(0)=+\infty$. For all $c \neq 0$, we have $f'(c)=(p-1)\abs{c}^{p-2}$. Thus applying Lemma~\ref{lem:genmvt}, we get the following corollary.

\begin{cor}\label{cor:genmvtlp}
Let $a < b$, both $a$ and $b$ being finite. Then, for any $p > 1$, there exists $c \in ]a,b[ \setminus \{0\}$ such that 
\[
\abs{b}^{p-2}b-\abs{a}^{p-2}a=(p-1)\abs{c}^{p-2}(b-a).
\]
\end{cor}

\section{On $L^p$ spaces inclusion}
Since $\O$ has finite Lebesgue measure, we have the classical inclusion $L^q(\Omega) \subset L^p(\Omega)$ for $1 \leq p \leq q < +\infty$. More precisely
\[
\norm{f}_{L^p(\O)} \leq \abs{\O}^{1/p-1/q}\norm{f}_{L^q(\O)} = \norm{f}_{L^q(\O)} \leq \norm{f}_{L^\infty(\O)} ,
\]
since $\abs{\O}=1$. We also have the following useful (reverse) bound whose proof is based on H\"older inequality.
\begin{lem}
For any $1 \leq q < p < +\infty$ we have
\[
\norm{f}_{L^p(\O)} \leq \norm{f}_{L^\infty(\O)}^{1-q/p}\norm{f}_{L^q(\O)}^{q/p} .
\]
In particular, for $p > 1$
\[
\norm{f}_{L^p(\O)} \leq \norm{f}_{L^\infty(\O)}^{1-1/p}\norm{f}_{L^1(\O)}^{1/p} .
\]
\label{inclusion}
\end{lem}



\paragraph{Acknowledgement.}
This work was supported by the ANR grant GRAPHSIP. JF was partly supported  by Institut Universitaire de France.

\bibliographystyle{abbrv}
\bibliography{biblio}

\end{document}